\documentclass[11pt,a4paper]{article}
\usepackage{graphicx} 
\usepackage{pifont}
\usepackage{bigints}
\usepackage[margin=1.70cm]{geometry}
\usepackage{authblk}

 \usepackage{setspace}
 \usepackage{lscape}
 \usepackage{array}
 \usepackage{tikz}
\usetikzlibrary{decorations.pathreplacing}
\usepackage{booktabs}
 \usepackage{multirow}
 \usepackage{amsmath}
 \usepackage{amsthm}
 \usepackage[utf8]{inputenc}
 \usepackage[pdfencoding=auto]{hyperref}
 \usepackage{longtable}
 \usepackage{amsfonts}
\usepackage{booktabs}
\usepackage{siunitx}
\newcolumntype{P}[1]{>{\raggedright\arraybackslash}p{#1}}
 \newtheorem{theorem}{Theorem}[section] 
 \newtheorem{proposition}{Proposition} [section]
 \newtheorem{definition}{Definition}[section]
 \newtheorem{example}{Example}[section]
  \newtheorem{cor}{Corollary}[section]
 \newtheorem{remark}{Remark}[section]
\numberwithin{equation}{section}

 \newtheorem{lemma}{Lemma}[section] 
 \usepackage{tabularx}
\usepackage{amssymb}
 \usepackage{caption}
 \usepackage{subcaption}
 \usepackage{array}
 \usepackage{xcolor}
 \usepackage{hyperref}
 \hypersetup{
    colorlinks=true,
    linkcolor=blue,
    filecolor=blue,      
    urlcolor=blue,
    citecolor=blue
}
 \usepackage[T1]{fontenc}
 \newcommand{\R}{\ensuremath{\mathbb{R}}}

\begin{document}

     \title{\bf A Generalized Notion of Completeness and Its Application}
    
         \date{}
      \author[1]{Himanshi Singh}
 \author[2]{Tanmay Sahoo}
 \author[1]{Nil Kamal Hazra}
\affil[1]{Department of Mathematics, Indian Institute of Technology Jodhpur, Rajasthan, India}
\affil[2]{Department of Mathematics, Indian Institute of Technology Palakkad,  Kerala, India}
   \maketitle

\begin{abstract}
{From the perspective of data reduction, the notions of minimal sufficient and complete statistics together play an important role in determining optimal statistics (estimators).  The classical notion of sufficiency and completeness are not adequate in many robust estimations that are based on different divergences. Recently, the notion of generalized sufficiency based on a generalized likelihood function was introduced in the literature. 
It is important to note that the concept of sufficiency alone does not necessarily produce optimal statistics (estimators). 
Thus, in line with the generalized sufficiency, we introduce a generalized notion of completeness with
respect to a generalized likelihood function. We then characterize the family of probability distributions that
possesses completeness with respect to the generalized likelihood function associated with the density power divergence (DPD). Moreover, we show that the family of distributions associated with the logarithmic density power divergence (LDPD) is not complete. Further, we
extend the Lehmann-Scheffé theorem and the Basu's theorem for the generalized likelihood estimation. Subsequently, we obtain the generalized uniformly minimum variance unbiased estimator (UMVUE)
for the $\mathcal{B^{(\alpha)}}$-family. Further, we derive an formula of the asymptotic expected deficiency (AED) that is used to compare the performance between the minimum density power divergence estimator (MDPDE) and the generalized UMVUE  for  $\mathcal{B^{(\alpha)}}$-family.  Finally, we provide an application of the developed results in stress-strength reliability model.
}
  \end{abstract}
\section{Introduction}\label{2.S1}
Let $Y_1^n := Y_1, \dots, Y_n$ be a random sample drawn from a population with the probability density function $f_\lambda \in \mathcal{P}$, where $\mathcal{P} = \{f_\lambda: \lambda \in \Lambda\}$ is a family of probability distributions parametrized by $\lambda$ (i.e., the functional form of $f_\lambda$ is known except the unknown parameter $\lambda$), and $\Lambda$ is an open subset of $\mathbb R^k$ and $\mathbb{R}$ is the set of real numbers. One of the key objectives in statistical inference is how one can use the information contained in $Y_1^n$ to draw meaningful inferences about $\lambda$. When the sample size $n$ is very large, the sample data $Y_1^n$ may be difficult to interpret directly. In these circumstances, statisticians often summarize the information in $Y_1^n$ using suitable statistics. Thus, any statistic can be viewed as a form of data reduction. Ideally, this process of data reduction through statistics should satisfy two key properties:
$(i)$ No important information about $\lambda$, present in $Y_1^n$,  should be lost;
$(ii)$ All irrelevant information in $Y_1^n$ that does not pertain to $\lambda$ should be discarded. In other words, a good statistic removes noise or extraneous features while preserving all important information about the parameter. Fisher~\cite{Sufficient} introduced the notion of sufficient statistic that preserves all information about the parameter $\lambda$ present in the sample $Y_1^n$. Further, the sufficient statistics for $\lambda$ can be obtained using the factorization theorem (see, \cite{HS, N}). Note that there can be more than one sufficient statistics for a given family. In fact, the whole sample $Y_1^n$ and order statistics are always sufficient statistics when the random sample is under consideration. One natural question arises which sufficient statistic should be preferred among the set of all sufficient statistics. 
The minimal sufficient statistic, introduced by Lehmann and Scheff\'e~\cite{LS}, is the one that provides the maximum data reduction while retaining all information about $\lambda$. A sufficient statistic $T$ is said to be minimal if $T$ is a function of all other sufficient statistics. Thus, the minimal sufficient statistic satisfies one of the key properties of the data reduction process. However, it may contain some noisy components, i.e., ancillary information (Lehmann and Scholz~\cite{LS1})). For example, let $Y_1^n = Y_1, \dots, Y_n$ be a random sample from $U(\lambda, \lambda+1)$. Then, $T(Y_1^n) = \left(Y_{(n)}- Y_{(1)}, \frac{Y_{(1)} + Y_{(n)}}{2}\right)$ is a minimal sufficient statistic for $\lambda$ but have ancillary component $Y_{(n)}- Y_{(1)}$. Thus, a minimal sufficient statistic does not always satistfy both the key properties of the data reduction process.
\\\hspace*{0.2 in}  In the context of point estimation, we often show interest to estimate some real-valued function $\tau(\lambda)$ for a given parametric family $\mathcal{P} = \{f_\lambda: \lambda \in \Lambda\}$. A sensible choice for an estimator $T$ is the one that minimizes the risk, which is usually quantified by a loss function. For example, under the squared error loss, a reasonable estimator minimizes the expected squared error, $E_{\lambda}[T -\tau(\lambda)]^2$, where $E_{\lambda}$ denotes the expectation with respect to the probability density function $f_\lambda$. Unfortunately, this estimator generally depends on the specific value of $\lambda$ at which the risk $E_{\lambda}[T -\tau(\lambda)]^2$ is minimized. If the true value is $\lambda_0$ then one could trivially choose $T= \tau(\lambda_0)$ to minimize the expected squared error. One strategy to get around this problem is to restrict the class of estimators by excluding estimators that lack particular desirable characteristics. One such desirable characteristic is the unbiasedness. An estimator is said to be unbiased estimator of $\tau(\lambda)$ if $E_{\lambda}[T(Y_1^n)] = \tau(\lambda)$, for all $\lambda \in \Lambda$. It frequently turns out that, among the set of all unbiased estimators, there is one that minimizes $E_{\lambda}[T -\tau(\lambda)]^2$, or equivalently, minimizes the variance of $T$, uniformly in $\lambda$ (i.e., uniformly minimum variance unbiased estimator (UMVUE)). However, finding the best estimator within the limited class of unbiased estimators is not always feasible. Then, an immediate question arises is $-$ under what conditions can such an estimator be obtained (see, Basu~\cite{B, B1})?
\\\hspace*{0.2 in} In order to accomplish the above mentioned query, Lehmann and Scheff\'e~\cite{LS} introduced the notion of completeness for a family of probability distributions. A family of probability distributions $\{f_\lambda: \lambda \in \Lambda\}$ is said to be complete (resp. boundedly complete) if, for every function (resp. bounded function) $\phi(\cdot)$ with $E_{\theta}[\phi(Y_1^n)]= 0$, for all $\lambda\in \Lambda$, we have $P_{\theta}\{\phi(Y_1^n) =0\} =1$, for all $\lambda\in \Lambda$. A sufficient statistic $T$ is said to be complete (resp. boundedly complete) if the family of probability distributions of $T$ is complete (resp. boundedly complete). 
It is worthwhile to note that the property of completeness plays a crucial role in hypothesis testing. In particular, it is relevant to similar regions, introduced by Neyman and Pearson~\cite{NP}, in the context of testing composite hypotheses. 
\\\hspace*{0.2 in} An important family of probability distributions which admits reduction by means of sufficient statistics is the exponential family, defined by probability densities of the form $f_\lambda(y) = \exp\left[h(y)+Z(\lambda)+{w(\lambda)}^T f(y)\right]$, $y \in S$, $\lambda \in \Lambda$. 
 It is worthwhile to note that, the entire sample is always a sufficient statistic and so, it is not helpful in data reduction. Thus,  in order to make a difference, sufficient statistics should be fixed in number, independent of sample size. Fisher, Darmois, Koopman, and Pitman~\cite{Fisher, Darmois, Koopman, Pitman} showed that, under certain regularity conditions, the existence of such a set of sufficient statistics is possible if and only if the underlying family is exponential. Thus, the exponential family has a fixed number of sufficient statistics (independent of sample size). Furthermore, if the exponential family has a full rank, the sufficient statistics are also complete. Recently, Gayen and Kumar~\cite{Generalized_sufficiency} further explored the notion of sufficient statistics for some generalized likelihood-based estimation. As mentioned in their paper, the classical notion of sufficiency is defined under the usual maximum likelihood (ML) estimation.
 However, a statistic that is sufficient in ML estimation may not be sufficient in other estimations. For example, $T(Y_1^n)=\overline Y$ is a sufficient statistic for the family of binomial distributions, Bin$(n, \lambda)$, when $n$ is known. But, this is no longer a sufficient statistic when the Cauchy-Schwarz
likelihood function is considered (for detailed explanation, see Gayen and Kumar~\cite{Generalized_sufficiency}). In this context, Gayen and Kumar~\cite{Generalized_sufficiency} introduced the notion of generalized sufficient statistic with respect to a generalized likelihood function. 
They have characterized families of probability distributions that possess a fixed number of sufficient statistics, regardless of sample size, with respect to the generalized likelihood functions associated with the density power divergence (DPD) and the logarithmic density power divergence (LDPD). Further, they extended the notion of minimal sufficiency with respect to a generalized likelihood function and derive generalized minimal sufficient statistics for the power-law family. Furthermore, they derived the generalized Rao-Blackwell theorem under a generalized likelihood estimation framework. 
It is important to note that the concept of sufficiency alone does not necessarily produce optimal statistics. The combination of sufficiency and completeness is often required to get optimal statistics, such as uniformly minimal variance unbiased estimator (UMVUE). 
Lehmann and Scheff\'e~(\cite{LS}, Theorem 5.1) stated that if $g(\lambda)$ is an estimable function and $T(Y_1^n)$ is a complete sufficient statistic for a given family $\mathcal P$, then $g(\lambda)$ has one and only one unbiased estimator that is a function of $T$, which is also a unique UMVUE of $g(\lambda)$. 
However, the existing notion of completeness will no longer be adequate to yield the desired results when working with generalized framework of sufficiency, especially one that goes beyond the usual likelihood-based method. This emphasizes the necessity of a generalized definition of completeness that is compatible with the notion of generalized sufficiency.
Further, from the brief discussion, one can conclude that completeness is a very important property for a family of probability distributions. Consequently, a few natural questions arise, as follows:
\begin{itemize}
\item How can one derive a generalized minimal sufficient statistic for a given family? 
\item When does a generalized minimal sufficient statistic contain no ancillary information?
\item How can one obtain a UMVUE for a given family with respect to a generalized likelihood function? 
\end{itemize}
To answer the questions mentioned above, in this paper, we generalize the notion of completeness with respect to
a generalized likelihood function. Then, we characterize the family of probability distributions that possesses completeness with respect to the generalized likelihood function associated with DPD. Further, we extend the Lehmann and Scheff\'e theorem for the generalized likelihood estimation. Subsequently, we derive the generalized UMVUE for the $\mathcal{B^{(\alpha)}}$-family. 

In a given estimation problem, there may have more than one estimators. In the literature, there are several measures that are used to decide which estimator can be be preferred over others. The asymptotic expected deficiency (AED), introduced by Hodges and Lehmann~\cite{HL}, is one of the popular measures that is often used. In this paper, we study the AED of the minimum density power divergence estimator (MDPDE) relative  to the generalized UMVUE for  $\mathcal{B^{(\alpha)}}$-family. 

In summary, the main objectives and the novelty of this paper are as follows:
\begin{enumerate}
    \item [$(i)$] We generalize the notion of completeness and ancillarity for generalized likelihood estimation; 
    \item [$(ii)$] We derive the generalized complete sufficient statistic for $\mathcal{B^{(\alpha)}}$-family. Further, we show that the underlying minimal sufficient statistic for $\mathcal{M^{(\alpha)}}$-family is not complete; 
    \item [$(iii)$] We extend the Lehmann-Scheffé theorem and the Basu's theorem for the generalized likelihood estimation. Moreover, we obtain the generalized UMVUE for $\mathcal{B^{(\alpha)}}$-family,
    \item [$(iv)$] We derive the AED of MDPDE relative to generalized UMVUE for $\mathcal{B^{(\alpha)}}$-family. Subsequently, we provide an application of this result in reliability stress-strength model.
\end{enumerate}
\hspace*{0.2 in}The rest of the paper is organized as follows. In Section \ref{2.S2}, we give some preliminaries and discuss some useful results. In Section \ref{2.S3}, we introduce the generalized notion of completeness and ancillarity, and study some associated results. 
In Section \ref{2.S5}, we study efficient estimators for some general families of distributions.
In Section \ref{arna0}, we compare the MDPDE with the generalized UMVUE  in terms of AED for $\mathcal{B^{(\alpha)}}$-family. Finally, some concluding remarks are given in Section~\ref{2.S6}.

\section{Preliminaries and some useful results}\label{2.S2}
In this section, we discuss some basic concepts that are used in subsequent sections. We first introduce some notation. We denote the set of real numbers and the set of positive real numbers by $\mathbb{R}$ and $\mathbb{R}^{+}$, respectively. We begin this section  by defining the notion of divergence, which measures the dissimilarity between two probability distributions (\cite{Basu_book}).
\begin{definition}\label{2.D2.1}
A divergence is a non-negative real-valued function $D$ such that, for any two probability density functions $f$ and $g$,
\begin{enumerate}
     \item [$(a)$] $D(g,f) \in \R$;
    \item [$(b)$] $D(g,f) \geq 0$;
    \item [$(c)$] $D(g,f) = 0$ if and only if $g=f$. $\hfill\Box$
    \end{enumerate} 
\end{definition}    
Various divergence measures exist in the literature (see \cite{ Entropy_Review, Divergence_review}). Some of the well-known ones include the Kullback-Leibler divergence (KLD), the Renyi divergence, the Tsallis divergence (power divergence (PD)), the density power divergence (DPD),  the relative $\alpha$-entropy (logarithmic density power divergence (LDPD)), the  S-divergence (SD), etc. (see \cite{DPD, Basu_book, LNRE, LDPD, KLD, Renyi, Tsallis}). In what follows, we give the definitions of some divergences that are used in this paper.
\begin{definition}\label{2.D2.2}
Let $f$ and $g$ be two probability density functions having the common support $S\subset \mathbb{R}$. Further, let $\alpha > 0$, $\alpha\neq 1$. Then, the
\begin{enumerate}
    \item [$(i)$] Kullback-Leibler divergence (KLD), denoted by $ \mathcal{D}_{KL} $, is given by 
\begin{equation}\label{2.E1}
  \mathcal{D}_{KL} (g,f)=  \int g(y) \log \frac{g(y)}{f(y)} dy; 
\end{equation}


\item [$(ii)$] Density power divergence $(DPD)$, denoted by $ \mathcal{B}_\alpha $, is given by
\begin{equation}\label{2.E4}
  \mathcal{B}_\alpha (g,f)= \frac{\alpha}{1-\alpha} \int g(y)f^{\alpha-1}(y)dy - \frac{1}{1-\alpha} \int g^\alpha(y) dy + \int f^{\alpha}(y) dy; 
\end{equation}
\item [$(iii)$] Logarithmic density power divergence (LDPD), denoted by $\mathcal{I}_\alpha $, is given by 
\begin{equation}\label{2.E7}
\mathcal{I}_\alpha (g,f)= \frac{\alpha}{1-\alpha} \log \int g(y)f^{\alpha-1}(y) dy - \frac{1}{1-\alpha} \log \int g^\alpha(y) dy + \log \int f^{\alpha}(y) dy. 
\end{equation}
\end{enumerate}
\end{definition}
\hspace*{0.02 in}Let $ Y_1^n := Y_1, \dots , Y_n$ be an independent and identically distributed (i.i.d.) random sample drawn according to a distribution having the probability density function (pdf) $g$, and let $y_1^n := y_1,\dots,y_n$ be the realized value of $Y_1^n$. Let $f_\lambda \in \mathcal{P}: = \{f_\lambda : \lambda \in \Lambda \subset \R^k\} $ be the model pdf with parameter $\lambda \in \Lambda$, which models the distribution $g$. Now, the parameter $\lambda$ can be estimated by minimizing the divergence $D(g, f_\lambda)$ with respect to $\lambda$. Then, the estimating equations corresponding to the divergences given in \eqref{2.E1}, \eqref{2.E4} and \eqref{2.E7} are given by
\begin{equation}\label{2.E8}
   \frac{1}{n} \sum_{i=1}^n u(y_j, \lambda) = 0,
\end{equation}

\begin{equation}\label{2.E9}
\frac{1}{n} \sum_{j=1} ^ {n} f_\lambda^{\alpha-1}(y_j)u(y_j,\lambda) = \int f_\lambda^\alpha(y)u(y;\lambda)  dy,
\end{equation}
and
\begin{equation}\label{2.E10}
\frac{\frac{1}{n} \sum_{j=1} ^ {n}f_\lambda^{\alpha-1}(y_j)u(y_j;\lambda)}{\frac{1}{n} \sum_{j=1} ^ {n} f_\lambda^{\alpha-1}(y_j)} = \frac{\int f_\lambda^\alpha(y)u(y;\lambda)  dy}{\int f_\lambda^\alpha(y)  dy},
\end{equation} 
respectively, where $u(y, \lambda)= \Delta \log  f_\lambda(y)$ (see \cite{DPD, LDPD, Windham}). 
\\\hspace*{0.2 in}The notion of sufficient statistics was first introduced by Fisher~\cite{Sufficient}. A sufficient statistic is a function of random sample and contains as much information about the unknown parameter of the underlying distribution as the sample has. An equivalent definition in terms of the usual likelihood function was given by Koopman~\cite{Koopman}. Recently, Gayen and Kumar~\cite{Generalized_sufficiency} introduced a generalized notion of sufficiency for estimation problems by considering the generalized likelihood function whose first-order optimality condition for maximization is the same as the estimating equation for minimization of the underlying divergence.
The generalized likelihood functions associated with the divergences given in  \eqref{2.E1}, \eqref{2.E4}, and \eqref{2.E7} 
are given by

\begin{equation}\label{2.E12}
 \mathcal{L}(y_1^n, \lambda) :=  \sum_{j=1}^n \log  f_\lambda (y_j),     
\end{equation}

\begin{equation}\label{2.E13}
 \mathcal{L}_{\mathcal{B}}^{(\alpha)}(y_1 ^n;\lambda):= \frac{1}{n}\sum_{j=1} ^ {n} \left[\frac{ \alpha f_\lambda^{\alpha-1}(y_j) - 1}{\alpha-1}\right] -  \int f_\lambda^\alpha(y)  dy, 
\end{equation}
and
\begin{equation}\label{2.E14}
 \mathcal{L}_{\mathcal{J}}^{(\alpha)}(y_1 ^n;\lambda) :=\frac{1}{\alpha-1}\log \left[\frac{1}{n}\sum_{j=1} ^ {n} f_\lambda^{\alpha-1}(y_j)\right] - \frac{1}{\alpha}\log \left[\int f_\lambda^\alpha(y)  dy\right],
\end{equation}
respectively. 

Below we give the definition of a deformed probability distribution, which is used throughout the paper (\cite{Generalized_sufficiency}).
\begin{definition}[Deformed probability distribution]\label{HTH}
Let $Y_1^n $ be a random sample drawn from a population having the probability density function $f_\lambda \in \mathcal{P}$ and let $y_1^n$ be a realized value of $Y_1^n $. Then, the deformed probability distribution associated with a generalized likelihood function $\mathcal{\mathcal{L_G}}$ is given by
\begin{equation}\label{2.E17}
    \widetilde{f}_\lambda(y_1^n) = \frac{\exp \left[\mathcal{L_G}(y_1^n; \lambda)\right]}{\int \exp[\mathcal{L_G}(r_1^n; \lambda)]dr_1^n},
\end{equation}
provided that the denominator is finite. Moreover, we denote this family of deformed distributions by $\widetilde{\mathcal{P}} = \{ \widetilde{f}_\lambda : \lambda \in \Lambda \}$, and write $\widetilde{P}_{\lambda}(\cdot)$, $\widetilde{E}_\lambda[\cdot]$ and $\widetilde{Var}_\lambda[\cdot]$ 
to represent the probability, the mean and the variance measures with respect to $\widetilde{f}_\lambda$, respectively. $\hfill\Box$
\end{definition}
In what follows, we give the definition of generalized sufficiency and discuss some associated results (\cite{Generalized_sufficiency}).
Before moving to the definition, we introduce some notation.
Let $\mathcal{P} = \{f_\lambda : \lambda \in \Lambda \}$ be a family of
probability distributions, and let $ Y_1^n := Y_1, \dots , Y_n$ be an i.i.d. random sample drawn according to some $f_\lambda\in \mathcal{P}$. Let $T := T(Y_1^n)$ be a statistic (which can be a vector valued function) and let $\mathcal{T} = \{ T(y_1^n): y_1^n \in S^n\subseteq \mathbb{R}^n\}$ be the set of values of $T$ for different realized samples. Suppose that the underlying problem is to estimate $\lambda$ by maximizing some generalized likelihood function $\mathcal{\mathcal{L_G}}$.
\begin{definition}[Generalized sufficient statistic]\label{GSS}
A statistic $T$ is said to be a generalized sufficient statistic for $\mathcal P$ if $[\mathcal{\mathcal{L_G}}(r_1^n; \lambda)-\mathcal{\mathcal{L_G}}(s_1^n; \lambda) ]$
is independent of $\lambda$ whenever $T(r_1^n)=T(s_1^n)$. 
\end{definition}
\begin{proposition}[Generalized factorization theorem]\label{GFT}
A statistic $T$ is said to be a generalized sufficient statistic for $\mathcal P$  if and only if there exist functions $p: \Lambda \times \mathcal{T} \rightarrow \R$ and $q: S^n \rightarrow \R$ such that 
\begin{equation}\label{2.E16}
    \mathcal{\mathcal{L_G}}(y_1^n; \lambda) = p(\lambda, T(y_1^n)) + q(y_1^n).
\end{equation} 
\end{proposition}
\begin{definition}[Generalized Fisher's definition]\label{GFD}
  A statistic $T$ is generalized sufficient for $\mathcal P$ if and only if 
\begin{equation}\label{2.E18}
\begin{aligned}
  \widetilde{f}_{\lambda_{y_1^n|t}}(y_1^n|t) & := \frac{\widetilde f_\lambda (y_1^n) \textbf{1}(T(y_1^n)=t)}{\widetilde{g}_\lambda(t)}  \\
  & = \begin{cases}
     \frac{\widetilde f_\lambda (y_1^n)}{\widetilde{g}_\lambda(t)} & \text{if $y_1^n\in C_t$}\\
    0 & \text{otherwise},
    \end{cases}
\end{aligned}
\end{equation}
is independent of $\lambda$, where $C_t = \{y_1^n : T(y_1^n)=t\}$ and
$\widetilde{g}_\lambda(t) := \int_{C_t} \widetilde f_\lambda (r_1 ^n)dr_1^n$. $\hfill\Box$
\end{definition}
 Next, we give the definition of generalized minimal sufficient statistic.

\begin{definition}[Generalized minimal sufficient statistic]\label{MSS}
A statistic $T$ is said to be a generalized minimal sufficient statistic for $\mathcal P$  if $T$ is a function of any other generalized sufficient statistics. Equivalently, $T$ is generalized minimal if, for any two i.i.d. samples $r_1^n$ and $s_1^n$, $T(r_1^n)=T(s_1^n)$ holds if and only if $[\mathcal{L_G}(r_1^n, \lambda) - \mathcal{L_G}(s_1^n, \lambda)]$ is independent of $\lambda$.  $\hfill\Box$ 
\end{definition}

\hspace*{0.02 in}Next, we discuss some general families of distributions. We first give the definition of the exponential family, which contains many popular distributions as special cases (\cite{exp}). This family is obtained as  the projection of KLD on a set  of probability distributions determined by linear constraints (\cite{Csiszar}).
\begin{definition}[Exponential family]\label{EF}
Let $ \{ f_\lambda : \lambda \in \Lambda \}$ be a family of probability distributions, where $\Lambda$ is an open subset of $\R^k$. Then, this family is said to form an exponential family, characterized by $w,h,f,\Lambda$ and $S$, if
\begin{equation}\label{2.E19}
  f_\lambda(y)=\begin{cases}
    \exp\left[h(y)+Z(\lambda)+{w(\lambda)}^T f(y)\right], & \text{for $y \in S$,}\\
    0, & \text{otherwise},
  \end{cases}
\end{equation}
where $S \subseteq \R$ is the support of $f_\lambda$; $e^{Z(\lambda)}$ is the normalizing constant; $w=[w_1, \dots ,w_d]^T$ and $w_i: \Lambda \rightarrow \R$ is differentiable, for $i= 1, \dots ,d$; $f=[f_1, \dots ,f_d]^T$ with $f_i:\R \rightarrow \R$, for $i= 1, \dots ,d$, and $h:\R\rightarrow \R$. Moreover, we denote this family by $\mathcal{E} $. $\hfill\Box$  
\end{definition}
Like the exponential family, there are some other generalized families that appeared as the projections of other divergences. We below give the definitions of ${\mathcal{B}^{(\alpha)}}$-family  and ${\mathcal{M}^{(\alpha)}}$-family.
\begin{definition}[$\mathcal{B}^{(\alpha)}$-Family]\label{B1F}
Let $ \{ f_\lambda : \lambda \in \Lambda \}$ be a family of probability distributions, where $\Lambda$ is an open subset of $\R^k$. Then, this family is said to form a $k$-parameter $ {\mathcal{B}^{(\alpha)}}$-family, characterized by $w,h,f,\Lambda$ and $S$,  if \begin{equation}\label{2.E20}
  f_\lambda(y)=\begin{cases}
    \left[h(y)+Z(\lambda)+ {w(\lambda)}^T f(y)\right]^\frac{1}{\alpha-1}, & \text{for $y \in S$,}\\
    0, & \text{otherwise},
  \end{cases}
\end{equation}
where $w$, $h$, $f$ and $S$ are the same as in Definition \ref{EF}, and $Z(\lambda) : \Lambda \rightarrow \R$ is the normalizing factor. Moreover, we denote this family by ${\mathcal{B}^{(\alpha)}}$.
\end{definition} 

\begin{remark}\label{2.R2.2}
The following observations can be made.
\begin{enumerate}
\item [$(a)$] If $f(y)= (1-\alpha)\widetilde{f}(y)$, $Z(\lambda) = (1-\alpha)Z'(\lambda) $ and $h(y)=q^{\alpha-1}(y)$ for some functions $\widetilde{f}$, $Z'$, $q$ ($> 0$), and $\alpha \rightarrow 1$, then $\mathcal{B}^{(\alpha)}$-family coincides with $\mathcal{E}$-family;
\item [$(b)$] The $ {\mathcal{B}^{(\alpha)}}$-family can be obtained through projection of DPD on a set of probability distributions determined by linear constraints (\cite{B-family}); 
 \item [$(c)$]  $\overline{f}^{d} = \Big[\overline{f_1},\dots,\overline{f}_{d}\Big]^T $ is a generalized sufficient statistic for the $k$-parameter $\mathcal{B}^{(\alpha)}$-family with respect to $\mathcal{L^{\alpha}_B}$, where  $\overline{f_i} = \frac{1}{n}\sum_{j=1}^{n}  f_i(y_j)$ for $i=1,\dots,d$ \cite{Generalized_sufficiency}. 

\end{enumerate} 
\end{remark}
\begin{definition}[$\mathcal{M}^{(\alpha)}$-family]\label{M1F}
Let $ \{ f_\lambda : \lambda \in \Lambda \}$ be a family of probability distributions, where $\Lambda$ is an open subset of $\R^k$. Then, this family is said to form a $k$-parameter $ {\mathcal{M}^{(\alpha)}}$-family, characterized by $w,h,f,\Lambda$ and $S$, if 
\begin{equation}\label{2.E21}
  f_\lambda(y)=\begin{cases}
    N(\lambda)\left[h(y)+ {w(\lambda)}^T f(y)\right]^\frac{1}{\alpha-1}, & \text{for $y \in S$}\\
    0, & \text{otherwise},
  \end{cases}
\end{equation}
where $w$, $h$, $f$ and $S$ are the same as in Definition \ref{EF}, and  $N(\lambda)^{-1} = \bigints \limits_S \left[h(y)+ {w(\lambda)}^T f(y)\right]^\frac{1}{\alpha-1}dy$. Moreover, we denote this family by ${\mathcal{M}^{(\alpha)}}$.
\end{definition}
\begin{remark}\label{2.R2.3} 
The following observations can be made.
\begin{enumerate}
    \item [$(a)$] If $f(y)= (1-\alpha)\widetilde{f}(y)$ and $h(y)=q^{\alpha-1}(y)$ for some function $\widetilde{f}$, $q$ ($> 0$) and $\alpha \rightarrow 1$, then $\mathcal{M}^{(\alpha)}$-family coincides with $\mathcal{E}$-family;
    \item [$(b)$] The $ {\mathcal{M}^{(\alpha)}}$-family can be obtained through projection of LDPD on a set of probability distributions determined by linear constraints (\cite{Forward_projection}, \cite{Reverse_projection});
    \item [$(c)$] Any member of $ {\mathcal{M}^{(\alpha)}}$-family  can be expressed as a member of $ {\mathcal{B}^{(\alpha)}}$-family and vice versa. However, this is not true for regular $ {\mathcal{B}^{(\alpha)}}$ and $ {\mathcal{M}^{(\alpha)}}$-families (\cite{Projection_theorems}).
    \item[$(d)$] ${\frac{\overline{f}}{\overline{h}}}^d = \left[\frac{\overline{f_1}}{\overline{h}},\dots. ,\frac{\overline{f_d}}{\overline{h}}\right]$ is a generalized sufficient statistic for the $k$-parameter $\mathcal{M}^{(\alpha)}$-family  with respect to $\mathcal{L}_{\mathcal{J}}^{(\alpha)}$, where $\overline{f_i} = \frac{1}{n}\sum_{j=1}^{n}  f_i(y_j)$ and $\overline{h} = \frac{1}{n}\sum_{j=1}^{n}  h(y_j)$ for $i=1,\dots,d$ (\cite{Generalized_sufficiency}).$\hfill\Box$\
\end{enumerate} 
\end{remark}
In a given estimation problem, the quality of an estimator is determined by how close it is to the true value of the parameter being estimated. The measure of closeness is determined through a loss function. Below we provide the definition of loss function and its corresponding risk.
\begin{definition}
Let $\mathcal{P} = \{f_\lambda : \lambda \in \Lambda \}$ be a family of
probability distributions, and let  $\mathcal{L_G}$  be the underlying generalized likelihood function. Further, let $\hat{\lambda}(Y_1^n)$ be an estimator of $\widetilde{\tau}(\lambda)$ such that $\widetilde{E}_\lambda(\hat{\lambda}(Y_1^n))=\widetilde{\tau}(\lambda)$. Then, $\hat{\lambda}(Y_1^n)$ is said to be a generalized unbiased estimator of $\widetilde{\tau}(\lambda)$. Moreover, $\widetilde{\tau}(\lambda)$ is said to be a generalized estimable function.
\end{definition}
\begin{definition}
Let $\mathcal{P} = \{f_\lambda : \lambda \in \Lambda \}$ be a family of
probability distributions, and let  $\mathcal{L_G}$  be the underlying generalized likelihood function.
Further, let $\widetilde{\tau}(\lambda)$ be a generalized estimable  function, and let it be estimated by $d \in \widetilde{\tau}(\Lambda)$. A function $L: \Lambda \times \widetilde{\tau}(\Lambda) \rightarrow \mathbb R^+$ is said to be a loss function, if 
$$L(\lambda, d) \geq 0, \text{ for all } (\lambda, d) \in \Lambda \times \widetilde{\tau}(\Lambda),$$and 
$$L[\lambda, \widetilde{\tau}(\lambda)] = 0, \text{ for all } \lambda \in \Lambda .$$
Then, the generalized risk of an estimator $\hat{\lambda}(Y_1^n)$ is defined as the generalized average loss, and is given by $$\widetilde R(\lambda, \hat{\lambda}) = \widetilde{E}_{\lambda}(L(\lambda, \hat{\lambda}(Y_1^n))).$$
\end{definition} 

\begin{remark}
Consider the squared error loss function given by
$$L(\lambda, d) = (d-\widetilde{\tau}(\lambda))^2, \quad (\lambda, d) \in \Lambda \times \widetilde{\tau}(\Lambda).$$
Then, for any estimator $\hat{\lambda}(Y_1^n)$ of a generalized estimable function $\widetilde{\tau}(\lambda)$, the risk reduces to generalized variance, i.e., $\widetilde R(\lambda, \hat{\lambda}) = \widetilde{Var}_\lambda(\hat{\lambda})$. $\hfill\Box$
\end{remark}
Below we provide the definition of the generalized uniformly minimum risk (resp. variance) unbiased estimator (UMRUE) (resp. UMVUE).
\begin{definition}
Let $\mathcal{P} = \{f_\lambda : \lambda \in \Lambda \}$ be a family of
probability distributions, and let  $\mathcal{L_G}$  be the underlying generalized likelihood function.
Further, let $\widetilde{\tau}(\lambda)$ be a generalized estimable function. An estimator $\hat{\lambda}(Y_1^n)$ is said to be generalized uniformly minimum risk (resp. variance) unbiased estimator of $\widetilde{\tau}(\lambda)$, if for any generalized unbiased estimator $\lambda^*(Y_1^n)$ of $\widetilde{\tau}(\lambda)$, we have
$$~~~~~~~~~~~~~~~~~~~~~~~~~~~~~~\widetilde R(\lambda, \hat{\lambda}) \;\leq \; \widetilde R(\lambda, \lambda^*) \;(\text{resp. } \widetilde{Var}_\lambda(\hat{\lambda}) \; \leq \;\widetilde{Var}_\lambda(\lambda^*)), \text{ for all } \lambda \in \Lambda.~~~~~~~~~~~~~~~~~~~~~~~~~~~~~~~~~~~\hfill\Box$$
\end{definition}
Sufficient statistics not only play an important role in data reduction but also are useful in constructing efficient estimators using the Rao-Blackwell theorem (\cite{Blackwell}, \cite{C.R. Rao}). Below we state the generalized Rao-Blackwell theorem which can be proved along the same lines as in Theorem 7.8 of \cite{Lehmann & Casella_book}.
\begin{theorem} \label{GRBT}
Let $\mathcal{P} = \{f_\lambda : \lambda \in \Lambda \}$ be a family of
probability distributions, and let  $\mathcal{L_G}$  be the underlying generalized likelihood function.
Further, let $T$ be a generalized sufficient statistic for $\mathcal{P}$. Let $\widetilde{\tau}(\lambda)$ be a generalized estimable function and $\hat{\lambda}$ be an  estimator of it such that $\widetilde{E}_\lambda[\hat{\lambda}(Y_1^n)]<\infty$. Further, let $L(\lambda, d)$ be a loss function which is  strictly convex in $d$. If 
$$\widetilde{R}(\lambda, \hat\lambda) =\widetilde{E}_\lambda[L(\lambda, \hat\lambda(Y_1^n)]<\infty, $$
and if 
$$\widetilde \varphi(T)=\widetilde{E}_\lambda[\hat{\lambda}(Y_1^n)|T],$$
then the generalized risk of the estimator $\widetilde \varphi(T)$ is less than $\widetilde{R}(\lambda, \hat\lambda)$ unless $\widetilde{P}_\lambda\left\{\hat{\lambda}(Y_1^n)=\widetilde \varphi(T)\right\}=1$.
\end{theorem}
\begin{remark}\label{HSH1}
    If we consider the loss function in the above theorem as $L(\lambda, d)=(d-\widetilde{\tau}(\lambda))^2$ and $\hat{\lambda}$ is a generalized unbiased estimator of $\widetilde{\tau}(\lambda)$, then $\widetilde{\varphi}(T)$ is also generalized unbiased estimator and have lesser variance than $\hat{\lambda}$ (see \cite{Generalized_sufficiency}, Theorem 18).
\end{remark}

\section{Generalized notion of completeness and ancillarity}\label{2.S3}
From the perspective of data reduction, the notions of minimal sufficient and complete statistics play an important role in determining the optimal statistics (estimators). Fisher's~\cite{Sufficient} concept of sufficiency is defined on the basis of the usual log-likelihood function, and the underlying estimation is the maximum likelihood (ML) estimation. However, in many practical scenarios, such as robust inference, it becomes necessary to move beyond ML-based methods. In such cases, the classical notion of sufficiency is not adequate. Recently, Gayen and Kumar~\cite{Generalized_sufficiency} introduced the notion of generalized sufficiency based on a generalized likelihood function. It is important to note that the concept of sufficiency alone does not necessarily produce optimal statistics. The combination of sufficiency and completeness is often required to get optimal statistics, such as UMVUE. However, the existing notion of completeness will no longer be adequate to yield the desired results when working with generalized framework of sufficiency, especially the one that goes beyond the usual likelihood-based method. This emphasizes the necessity of a generalized definition of completeness that is compatible with the notion of generalized sufficiency.

Thus, in line with the notion of generalized sufficiency, we introduce a generalized  notion of completeness and ancillarity. 
Before giving the definition of a generalized complete  statistic, we below give the definition of a usual complete statistic.
\begin{definition}\label{complete}
Let $\mathcal{P} = \{ f_\lambda : \lambda \in \Lambda \}$ be a family of probability distributions. A statistic $T$ is said to be complete  for $\mathcal{P} $ in the usual sense (i.e., with respect to the usual likelihood function $\mathcal{L}$) if, for any function~$h$,
   $$ E_{\lambda}[h(T)] = 0, \quad \text{for all} \quad \lambda \in \Lambda,$$
 implies 
 $$~~~~~~~~~~~~~~~~~~~~~~~~~~~~~~~~~~~~~~~~~~~~{P}_{\lambda}\{h(T) = 0\} = 1,  \quad \text{for all} \quad \lambda \in \Lambda.~~~~~~~~~~~~~~~~~~~~~~~~~~~~~~~~~~~~~~\Box $$
\end{definition}
\hspace*{0.02 in}From the Fisher-Darmois-Koopman-Pitman theorem, we know that, for the ML estimation,  the exponential family is the only one that has a fixed number of sufficient statistics, independent of sample size. Consequently, by using the notion of sufficiency, completeness and ancillarity, we obtain the efficient estimator for the family of exponential distributions. However, this cannot be done for other generalized families of distributions, namely, $\mathcal{B^{(\alpha)}}$-family and $\mathcal{M^{(\alpha)}}$-family. Recently, Gayen and Kumar~\cite{Generalized_sufficiency} derived the generalized Fisher-Darmois-Koopman-Pitman theorem for some generalized likelihood estimations. Consequently, it triggers whether the UMVUE for the aforementioned generalized families of distributions can be derived similar to that of the exponential family. To address this problem, we below introduce the notion of generalized complete statistic and generalized ancillary statistic. 
\\\hspace*{0.2 in}For a given generalized likelihood estimation problem, a generalized complete statistic plays the same role as the usual complete statistic does for the ML estimation. A generalized complete statistic is defined based on the deformed distribution associated with some generalized likelihood function $\mathcal{L_G}$. Below we give the definition of a generalized complete statistic.
\begin{definition}\label{2.D3.1} 
Let $\mathcal{P} = \{f_\lambda : \lambda \in \Lambda \}$ be a family of
probability distributions, and let  $\mathcal{L_G}$  be the underlying generalized likelihood function.
A statistic $T$ is said to be generalized complete  for $\mathcal{P} $ if, for any function~$h$, 
 $$\widetilde{E}_{\lambda} [h(T)] = \int \frac{ h(T(y_1^n))\exp \left[\mathcal{L_G}(y_1^n; \lambda)\right]}{\int \exp[\mathcal{L_G}(r_1^n; \lambda)]\; dr_1^n}\; dy_1^n = 0,  \quad \text{for all} \quad \lambda \in \Lambda,$$
 implies 
$$~~~~~~~~~~ \widetilde{P}_{\lambda}\{h(T) = 0\} =\int\limits_{\{y_1^n:h(T(y_1^n))=0\}} \frac{ \exp \left[\mathcal{L_G}(y_1^n; \lambda)\right]}{\int \exp[\mathcal{L_G}(r_1^n; \lambda)]\; dr_1^n}\; dy_1^n= 1,  \quad \text{for all} \quad \lambda \in \Lambda.~~~~~~~~~~~~~~~~$$  
Moreover, a statistic $T$ is said to be generalized complete sufficient statistic for $\mathcal{P} $ if it is both generalized sufficient and generalized complete. $\hfill\Box$
\end{definition}
 Next, we give the definition of a generalized ancillary statistic for the generalized likelihood estimation.  An ancillary statistic plays the opposite role to that of a  sufficient statistic. A  sufficient statistic for a given family of distributions contains all information about the unknown parameter(s), whereas an ancillary statistic does not contain any information about the parameter(s), i.e., its distribution is free of parameter(s).
\begin{definition}\label{2.D3.2}
Let $\mathcal{P} = \{f_\lambda : \lambda \in \Lambda \}$ be a family of
probability distributions, and let  $\mathcal{L_G}$  be the underlying generalized likelihood function.
  A statistic $A(Y_1^n)$ is said to be generalized ancillary for $\mathcal{P} $ if 
  $$~~~~~~~~~~~~~~~~~~~~\widetilde{P}_\lambda\{A(Y_1^n)\leq a\} = \bigintsss \limits_{y_1^n: A(y_1^n)\leq a} \frac{ \exp \left[\mathcal{L_G}(y_1^n; \lambda)\right]}{\int \exp[\mathcal{L_G}(r_1^n; \lambda)]\; dr_1^n} \;dy_1^n , \text{ is independent of } \lambda.~~~~~~~~~~~~~~~~~~~~~~~~\hfill\Box$$
\end{definition}
\hspace*{0.02 in}Both sufficiency and completeness are important in constructing the UMVUE for a given family of probability distributions. By considering the generalized Rao-Blackwell theorem and the generalized notion of completeness, we below give the generalized Lehmann-Scheffé theorem.  
\begin{theorem}\label{2.T3.1}
 Let $\mathcal{P} = \{f_\lambda : \lambda \in \Lambda \}$ be a family of
probability distributions, and let  $\mathcal{L_G}$  be the underlying generalized likelihood function. Further, let $T$ be a generalized complete sufficient statistic for $\mathcal{P}$. Then, every generalized estimable function $\widetilde{\tau}(\lambda)$ has one and only one generalized unbiased estimator that is of the form $h(T)$, for a fixed function $h$. Moreover, $h(T)$ is a unique generalized UMVUE of $\widetilde{\tau}(\lambda)$.
 
\end{theorem}

\noindent\textbf{Proof:} Since $\widetilde{\tau}(\lambda)$  is a generalized estimable function, there exists an estimator $\hat\lambda$ such that $\widetilde{E}_{\lambda} [\hat \lambda]=\widetilde{\tau}(\lambda) $. Then, $\widetilde{\varphi}(T) := \widetilde{E}_{\lambda}[\hat\lambda | T=t]$ is a generalized unbiased estimator for $\widetilde{\tau}(\lambda)$. If possible, let there be another generalized unbiased estimator $k(T)$ that is a function of $T$.
\noindent Then, we get
\begin{equation*}
     \widetilde{E}_{\lambda}[\widetilde{\varphi}(T)-k(T)]= 0, \quad \text{for all} \quad \lambda \in \Lambda.
\end{equation*}
Since T is a generalized complete sufficient statistic,  the above statement implies that 
\begin{equation*}
     \widetilde{P}_{\lambda}\{\widetilde{\varphi}(T)) - k(T) = 0\} = 1,  \quad \text{for all} \quad \lambda \in \Lambda ,  
\end{equation*}
or equivalently,
\begin{equation*}
\widetilde{ P}_{\lambda}\{\widetilde{\varphi}(T) = k(T) \} = 1,  \quad \text{for all} \quad \lambda \in \Lambda.   
\end{equation*}
Thus, $\widetilde{\varphi}(T)$ and $ k(T)$ differ only on a set of measure zero. Hence there exists one and only one generalized unbiased estimator, which is a function of $T$, of $\widetilde{\tau}(\lambda)$.
Now, by Theorem~\ref{GRBT}, we get 
\begin{equation*}
    \widetilde{Var}_{\lambda}(\widetilde{\varphi}(T) ) \leq \widetilde{Var}_{\lambda}(\hat \lambda), \quad \text{for all} \quad \lambda \in \Lambda.
\end{equation*}
As $\hat\lambda$ is arbitrary, $\widetilde{\varphi}(T)$ is a generalized UMVUE of $\widetilde{\tau}(\lambda)$, and it is unique.
~$\hfill\Box$

Below we state the same result as in the above theorem in terms of general loss function. The proof follows along the same lines.
\begin{theorem}
Let $\mathcal{P} = \{f_\lambda : \lambda \in \Lambda \}$ be a family of
probability distributions, and let  $\mathcal{L_G}$  be the underlying generalized likelihood function. Further, let $T$ be a generalized complete sufficient statistic for $\mathcal{P}$.
\begin{enumerate}
    \item [$(a)$]For every generalized estimable function $\widetilde{\tau}(\lambda)$, there exists a generalized UMRUE for any loss function $L(\lambda,d)$ which is convex in~$d$;
    \item [$(b)$] There is a unique generalized UMRUE that is a function of $T$, provided its risk is finite and $L$ is strictly convex in $d$. $\hfill\Box$
\end{enumerate}
\end{theorem}

One of the important results in estimation theory is the Basu's theorem, which states that any ancillary statistic is independent of the complete sufficient statistic with respect to the usual likelihood function (\cite{D_Basu}). In the following theorem, we show that this result indeed holds for the generalized likelihood estimation.
\begin{theorem}\label{2.T3.2}
Let $\mathcal{P} = \{f_\lambda : \lambda \in \Lambda \}$ be a family of
probability distributions, and let  $\mathcal{L_G}$  be the underlying generalized likelihood function.
Further, let $T$ be a generalized complete sufficient statistic for $\mathcal{P}$. Then, any generalized ancillary statistic $A$ for $\mathcal{P}$ is independent of $T$.  
\end{theorem}
\noindent\textbf{Proof:} Since $A$ is a generalized ancillary statistic, we have that
$$\widetilde{P}_\lambda\{A(Y_1^n)\leq a\} = \bigintsss \limits_{y_1^n: A(y_1^n)\leq a} \frac{ \exp \left[\mathcal{L_G}(y_1^n; \lambda)\right]}{\int \exp[\mathcal{L_G}(r_1^n; \lambda)]\; dr_1^n} \;dy_1^n , $$
is independent of $\lambda$, for all $a$.
Let $$h_a (t) = \widetilde{P}_\lambda \{A(Y_1^n)\leq a | T(Y_1^n) = t\}.$$
Now, by using the law of total probability, we get
$$ \widetilde{E}_\lambda[h_a(T(Y_1^n))] = \widetilde{P}_\lambda \{A(Y_1^n)\leq a\},$$
which implies
$$ \widetilde{E}_\lambda\left[h_a(T(Y_1^n)) - \widetilde{P}_\lambda \{A(Y_1^n)\leq a\}\right] = 0.$$
Since $T$ is a generalized complete sufficient statistic, we get from the above statement that
$$\widetilde{P}_\lambda \left\{ h_a(T(Y_1^n)) - \widetilde{P}_\lambda \left\{A(Y_1^n)\leq a\right\} = 0 \right\}=1,$$
or equivalently,
$$\widetilde{P}_\lambda \left\{A(Y_1^n)\leq a | T(Y_1^n) = t\right\} = \widetilde{P}_\lambda \left\{A(Y_1^n)\leq a\right\}$$
with probability $1$. Hence $A$ and $T$ are independent. $\hfill\Box$
\\\hspace*{0.2 in}In the next theorem, we establish a relationship between a complete sufficient statistic and a minimal sufficient statistic for the generalized likelihood estimation.
\begin{theorem}\label{2.T3.3}
Let $\mathcal{P} = \{f_\lambda : \lambda \in \Lambda \}$ be a family of
probability distributions, and let  $\mathcal{L_G}$  be the underlying generalized likelihood function.
If a generalized minimal sufficient statistic exists for $\mathcal{P}$, then any generalized complete sufficient statistic is also a generalized minimal sufficient statistic for $\mathcal{P}$.
\end{theorem} 

\noindent\textbf{Proof:} 
 Let $T_0$ be a generalized minimal sufficient statistic with $ \widetilde{E}_{\lambda}[T_0^2] < \infty$, and let $T$ be a generalized complete sufficient statistic. Then, by the definition of the generalized minimal sufficient statistic, we get $T_0 = h(T)$ with probability $1$, for some function $h$.
Further, let $ \psi(T_0): = \widetilde{E}_{\lambda}[T|T_0] $. Then, 
we get $ \widetilde{E}_{\lambda}[T - \psi(h(T))] =0$, for all $\lambda \in \Lambda$, which further implies that $T = \psi(h(T)) = \psi(T_0)$ with probability $1$.  
 Hence $T$ is a generalized minimal sufficient statistic.
\section{Efficient estimators for some general families of distributions}\label{2.S5}
In this section, we first discuss the necessary and sufficient conditions for an estimator to be generalized UMVUE for a given family of distributions with respect to a generalized likelihood function.
Then, we obtain the generalized UMVUE for $\mathcal{B}^{(\alpha)}$-family with respect to $\mathcal{L}_{\mathcal{B}}^{(\alpha)}$. Further, we verify whether the minimal sufficient statistic ${\frac{\overline{f}}{\overline{h}}}^d$ for $\mathcal{M}^{(\alpha)}$-family is the generalized UMVUE for the given family with respect to  
$\mathcal{L}_{\mathcal{J}}^{(\alpha)}$. 
\\\hspace*{0.2 in}We begin this section with the following theorem wherein we give the necessary and sufficient conditions for an estimator to be generalized UMVUE.  
\begin{theorem}\label{2.T5.6}
Let $\mathcal{P} = \{f_\lambda : \lambda \in \Lambda \}$ be a family of
probability distributions, and let  $\mathcal{L_G}$  be the underlying generalized likelihood function.
Let $\widetilde{W} = \{ T : \widetilde{E}_{\lambda} [T] = \widetilde{\tau}(\lambda), \widetilde{E_{\lambda}} [T^2] < \infty \;\text{ for all } \;\lambda \in \Lambda\}$ and $\widetilde{W}_0 = \{ w :  \widetilde{E}_{\lambda}[w] = 0, \widetilde{E}_{\lambda}[w^2] < \infty \text{ for all } \lambda \in \Lambda\}$. 
Then, $T_0 \in \widetilde{W}$ is the generalized UMVUE for $\lambda$ if and only if 

\begin{equation}\label{2.E52}
   \widetilde{E}_{\lambda}[wT_0]=0, \text{ for all }  \lambda\in \Lambda \text{ and } \text{ for all }  w \in \widetilde{W}_0.  
\end{equation}   
\end{theorem}

\noindent\textbf{Proof:} Let $T_0 \in \widetilde{W}$ be the generalized UMVUE. If possible, let   $\widetilde{E}_{\lambda_0}[w_0T_0] \neq 0 $ for some $\lambda_0$ and for some $ w_0 \in \widetilde{W}$. Note that, if $\widetilde{E}_{\lambda_0} [w_0^2] = 0$ then $\widetilde{P}_{\lambda_0}\{w_0=0\}=1$ and consequently, $ \widetilde{E}_{\lambda_0}[w_0T_0]=0 $. Thus, we assume that $\widetilde{E}_{\lambda_0} [w_0^2] > 0$.  Consider a statistic $T_0+\lambda w_0$ for some $\lambda\in R$. Clearly $T_0+\lambda w_0 \in \widetilde{W}$. Let us choose $\lambda = -{\widetilde{E}_{\lambda_0}[w_0T_0]}/{\widetilde{E}_{\lambda_0} [w_0^2]}$. Then,
\begin{equation*}
\widetilde{E}_{\lambda_0}[T_0+\lambda w_0]^2 = \widetilde{E}_{\lambda_0}[T_0^2]-\frac{(\widetilde{E}_{\lambda_0}[w_0T_0])^2}{\widetilde{E}_{\lambda_0}[w_0^2]} < \widetilde{E}_{\lambda_0}[T_0^2],
\end{equation*}
which implies that
\begin{equation*}
    \widetilde{Var}_{\lambda_0}(T_0+\lambda w_0)< \widetilde{Var}_{\lambda_0}(T_0).
\end{equation*}
Note that this is a contradiction.
Conversely, let $\widetilde{E}_{\lambda}[wT_0]=0 $, for some $T_0 \in \widetilde{W}$, $\text{for all } \lambda \in \Lambda$ and $\text{for all } w \in \widetilde{W}_0$, and let $T \in \widetilde{W}$. Then $T_{0}-T \in \widetilde{W_0}$ and 
$\widetilde{E}_{\lambda}[T_0(T_0 - T)]=0 $. By using Cauchy-Schwarz inequality, we get
\begin{equation*}
 \widetilde{E}_{\lambda}[T_0^2]=\widetilde{E}_{\lambda}[T_0 T] \leq (\widetilde{E}_{\lambda}T_0 ^2)^{1/2}(\widetilde{E}_{\lambda}T^2)^{1/2}.
\end{equation*}
If $\widetilde{E}_{\lambda}[T_0^2]=0$ implies $\widetilde{P}_{\lambda}\{T_0 = 0\} = 1$, there is nothing to prove. Otherwise
\begin{equation*}
   \widetilde{Var}_{\lambda}(T_0) \leq \widetilde{Var}_{\lambda}(T).
\end{equation*}
Since $T$ is arbitrary,  $T_0$ is the generalized UMVUE. Hence the result.
\subsection{$\mathcal{B^{(\alpha)}}$-family}
\hspace*{0.02 in}
We begin this subsection with the following theorem, wherein we obtain a generalized complete sufficient statistic for ${\mathcal{B}^{(\alpha)}}$-family.
\begin{theorem}\label{2.T5.4}
Consider the ${\mathcal{B}^{(\alpha)}}$-family, as in Definition~\ref{B1F}, with common support $S\subset \mathbb{R}$, and let $\mathcal{L_{B}^{(\alpha)}}$ be the underlying generalized likelihood function. Suppose that the range of $w(\lambda)$ contains a $d$-dimensional rectangle. Then, the generalized sufficient statistic $\overline{f}^{d} = \Big[\overline{f_1},\dots,\overline{f}_{d}\Big]^T $, as given in Remark~\ref{2.R2.2}(c), is generalized complete for the ${\mathcal{B}^{(\alpha)}}$-family. 
\end{theorem}

\noindent\textbf{Proof:} 
The deformed probability distribution associated with $\mathcal{L_{B}^{(\alpha)}}$ is given by
\begin{equation*}
\begin{aligned}
\widetilde{f}_\lambda(y_1^n) & = \frac{\exp\left[\mathcal{L^\alpha_B}(y_1^n;\lambda)\right]}{\bigintss \exp\left[\mathcal{L^\alpha_B}(r_1^n;\lambda)\right] dr_1^n} \\
& = \frac{\exp\left[\frac{\alpha}{n(\alpha-1)}\sum_{j=1}^n f_\lambda(y_j)^{\alpha-1}\right]}{\bigintss \exp\left[\frac{\alpha}{n(\alpha-1)}\sum_{j=1}^n f_\lambda(r_j)^{\alpha-1}\right] dr_1^n}.
\end{aligned}
\end{equation*}
After putting the value of $f_\lambda$ in the above expression, we get
\begin{equation}\label{2.E46}
\widetilde{f}_\lambda(y_1^n) = N(\lambda)\exp\left[ \frac{\alpha}{\alpha-1} \overline{h}(y_1^n)+\frac{\alpha}{\alpha-1} \sum_{i=1}^d w_i(\lambda)\overline{f}_i (y_1^n)\right],
\end{equation}
where $\overline{h}(y_1^n)=\frac{1}{n}\sum_{j=1}^n h(y_j)$, $ \overline{f_i}(y_1^n) = \frac{1}{n}\sum_{j=1}^n f_i(y_j)$ and $N(\lambda)^{-1} = \bigintsss \exp\Big[ \frac{\alpha}{\alpha-1}\overline{h}+\frac{\alpha}{\alpha-1}w^T(\lambda)\overline{f}^d\Big] d r_1^n $.
Further, the above expression can equivalently be written as
\begin{equation}\label{2.E47}
\widetilde{f}_\lambda(y_1^n) = N(\lambda)p(y_1^n) \exp\left[ \sum_{i=1}^d w^*_i(\lambda)\overline{f}_i (y_1^n)\right],
\end{equation}
where $w^*_i(\lambda)= \frac{\alpha}{\alpha-1}w_i(\lambda)$ and  $p(y_1^n) = \exp\left[ \frac{\alpha}{\alpha-1} \overline{h}(y_1^n)\right]$.
Next, by applying the generalized factorization theorem (Proposition~\ref{GFT}), we get that $\overline{f}^{d}=\left[\overline{f_1},\dots, \overline{f_d}\right]^T$ is a generalized sufficient statistic for $\lambda \in \Lambda$.  Further, by absorbing the factor $p(y_1^n)$ into a measure $\nu$, the  expression given in \eqref{2.E47} can be written as
\begin{equation*}\label{2.E49}
    d\widetilde{F}_\lambda(y_1^n) = N(\lambda)\exp\left[ \sum_{i=1}^d w^*_i(\lambda)\overline{f}_i (y_1^n)\right] d\nu(y_1^n). 
\end{equation*}
Since the range of $w(\lambda)$ contains a $d$-dimensional rectangle, the range of $w^*(\lambda) = \left(w^*_1(\lambda),\dots,w^*_{d}(\lambda)\right)$ also contains a $d$-dimensional rectangle. Consequently, by Theorem 1 in Section 4.3 of Lehmann~\cite{Lehmann_book},  we get that $\overline{f}^{d} = \Big[\overline{f_1},\dots,\overline{f}_{d}\Big]^T $ is a generalized complete sufficient statistic.$\hfill\Box$
\\\hspace*{0.2 in}The following corollary follows from Theorems \ref{2.T3.3} and \ref{2.T5.4}.
\begin{cor}
Under the assumption of Theorem \ref{2.T5.4}, 
$ \overline{f}^{d} = \Big[\overline{f_1},\dots,\overline{f}_{d}\Big]^T$ is a generalized minimal sufficient statistic for $\mathcal{B}^{(\alpha)}$-family.  $\hfill\Box$
\end{cor}
\hspace*{0.2 in}In the next theorem, we obtain the generalized UMVUE for $\mathcal{B}^{(\alpha)}$-family. The proof follows from Theorems~\ref{2.T3.1} and \ref{2.T5.4}.
\begin{theorem}\label{2.T5.5}
Consider the ${\mathcal{B}^{(\alpha)}}$-family, as in Definition~\ref{B1F}, with common support $S\subset \mathbb{R}$,  and let $\mathcal{L_{B}^{(\alpha)}}$ be the underlying generalized likelihood function. Suppose that the range of $w(\lambda)$ contains a $d$-dimensional rectangle. Then, for any function $g$, $g(\overline{f}^{d})$ is the generalized UMVUE for $\widetilde{\tau}(\lambda) = \widetilde{E_\lambda}\left[g(\overline{f}^{d})\right]$.  $\hfill\Box$
\end{theorem}
\hspace*{0.2 in}Below we give an example wherein we obtain the generalized UMVUE for the family of Student distributions with respect to $\mathcal{L}_{\mathcal{B}}^{(\alpha)}$.
\begin{example}\label{2.Ex5.1}
 Consider the Student distribution with the probability density function (pdf) given by
\begin{equation*}
  f_\lambda(y) = \frac{\Gamma(\frac{\nu+1}{2})}{\Gamma(\frac{\nu}{2})\sqrt{\pi \nu \sigma^2}} \Big[1+\frac{1}{\nu} \Big(\frac{y-\mu}{\sigma}\Big)^2\Big]^{-\frac{\nu+1}{2}}, 
  \end{equation*}
  where $\nu$ ($ > 2$) is the degrees of freedom and $\lambda = (\mu, \sigma^2) \in \R \times \R^+$ is the parameter set. By rearranging the terms, the above 
pdf can equivalently be written as
\begin{equation*}
   f_\lambda(y) =  \left[Z_\nu(\lambda)+  \frac{N_\nu}{\nu\sigma^{2\nu}}y^2 - \frac{2 \mu N_\nu}{\nu \sigma^{2\nu} }y \right]^{-\frac{\nu+1}{2}},   
\end{equation*}
where $ N_\nu=\left(\frac{\Gamma(\frac{\nu+1}{2})}{\Gamma(\frac{\nu}{2})\sqrt{\pi \nu }}\right)^{-\frac{2}{\nu+1}}$ and $Z_\nu(\lambda)=\frac{N_\nu}{(\sigma)^{-\frac{2}{\nu+1}}}\left[1+\frac{\mu^2}{\nu\sigma^2}\right]$. Further, this can equivalently be written as
\begin{equation}
 f_\lambda(y) =  \left[Z_{l_\alpha}(\lambda)+ \frac{N_{l_\alpha}}{\nu \sigma^{2l_\alpha}}y^2 -  \frac{2 \mu N_{l_\alpha}}{ \nu\sigma^{2l_\alpha} }y\right]^{\frac{1}{\alpha-1}},    
\end{equation}
where  $\alpha = 1 -\frac{2}{\nu+1}$, $ l_{\alpha} = {\nu} =\frac{1+\alpha}{1-\alpha}>2$. By considering  $h(y)=0$, $f(y) = [y^2, y]^T$ and $w(\lambda) = \left[\frac{N_{l_{\alpha}}}{ \nu\sigma^{2l_{\alpha}}}, \frac{-2\mu N_{l_{\alpha}}}{\nu \sigma^{2 l_{\alpha}}}\right]^T$, we note that $f_\lambda$ belongs to the $\mathcal{B}^{(\alpha)}$-family. Clearly, the range of $w(\lambda)$ contains a two-dimensional rectangle. Thus, by Theorem \ref{2.T5.4}, we get that $T = \left[\frac{1}{n} \sum_{j=1}^n y_j^2, \frac{1}{n} \sum_{j=1}^n y_j \right]^T$ is a complete sufficient statistic for $\lambda = (\mu, \sigma^2)$ with respect to $\mathcal{L}_{\mathcal{B}}^{(\alpha)}$. Consequently,  by Theorem \ref{2.T5.5}, we get that $g(T)$, for any function $g$, is the UMVUE for $\widetilde{\tau}(\lambda)=\widetilde{E}_\lambda[g(T)]$ with respect to $\mathcal{L}_{\mathcal{B}}^{(\alpha)}$.  
\end{example}

\subsection{$\mathcal{M}^{(\alpha)}$-family}
It is mentioned in Remark 3 of Gayen and Kumar~\cite{Generalized_sufficiency} that ${\frac{\overline{f}}{\overline{h}}}^d = \left[\frac{\overline{f_1}}{\overline{h}},\dots. ,\frac{\overline{f_d}}{\overline{h}}\right]$ is a generalized sufficient statistic for the $k$-parameter $\mathcal{M}^{(\alpha)}$-family (as in Definition~\ref{M1F}) with respect to $\mathcal{L}_{\mathcal{J}}^{(\alpha)}$, where $\overline{f_i} = \frac{1}{n}\sum_{j=1}^{n}  f_i(y_j)$ and $\overline{h} = \frac{1}{n}\sum_{j=1}^{n} h(y_j)$ for $i=1,\dots,d$. Moreover, it was also shown in Theorem 12 of Gayen and Kumar~\cite{Generalized_sufficiency} that this is indeed a minimal statistic for the regular $\mathcal{M}^{(\alpha)}$-family. In line with the result given in Theorem~\ref{2.T5.4}, one may be interested to know whether ${\frac{\overline{f}}{\overline{h}}}^d$ is a generalized complete statistic for $\mathcal{M^{(\alpha)}}$-family with respect to $ \mathcal{L^{(\alpha)}_{J}}$. Below we give a counterexample which shows that this is not the case.

\begin{example}\label{2.Ex5.2}
Consider the family of Bernoulli distributions with the pdf given by
\begin{equation*}
    f_\lambda(y) = (1-\lambda)\Big[1+y\frac{2\lambda-1}{1-\lambda}\Big], \quad y \in \{0, 1\} \quad \text{and} \quad \lambda \in (0,1).
\end{equation*}
Note that  $f_\lambda\in \mathcal{M}^{(\alpha)}$-family when $\alpha = 2$, $h(y)= 1$, $Z(\lambda) = 1-\lambda$, $w(\lambda) = \frac{2\lambda-1}{1-\lambda} $ and $f(y) = y$. Let $Y_1,Y_2,Y_3$ be a random sample drawn from $f_\lambda$. Then, the generalized likelihood function associated with the LDPD is given by
\begin{eqnarray}\label{2.E51}
\mathcal{L^{(\alpha)}_{J}}(y_1,y_2,y_3;\lambda) &=& \log \Bigg[\frac{\frac{1}{3}\sum_{j=1} ^3 f_\lambda (y_j)}{(E_\lambda[f_\lambda(Y)])^{1/2}}\Bigg] \nonumber \\
&=& \log \left[(1-\lambda) + (2\lambda-1)\frac{1}{3}\sum_{j=1} ^3 y_j \right]- \log\left[E_\lambda[f_\lambda(Y)]\right]^{1/2}.
\end{eqnarray}
Consequently, by using \eqref{2.E17}, the deformed pdf  associated with $\mathcal{L^{(\alpha)}_{J}}$ is given by
\begin{equation*}
 \widetilde{f}_{\lambda}(y_{1},y_2,y_3) = \frac{1+\overline{y}\frac{2\lambda - 1}{1-\lambda}}{\sum\limits_{r_1,r_2,r_3}\left(1+\overline{r}\frac{2\lambda - 1}{1-\lambda} \right) }= Z(\lambda) \Bigg[1+\overline{y}\frac{2\lambda - 1}{1-\lambda}\Bigg],  
\end{equation*}
where $\overline{y}=\frac{1}{3}\sum_{j=1}^3 y_j$ and $Z^{-1}(\lambda)=\sum\limits_{r_1,r_2,r_3}\left(1+\overline{r}\frac{2\lambda - 1}{1-\lambda}\right)$. Note that $\widetilde{f}_{\lambda}$ also belongs to $\mathcal{M^{(\alpha)}}$-family. Now, by  Theorem 12 of Gayen and Kumar~\cite{Generalized_sufficiency}, we get that ${\frac{\overline{f}}{\overline{h}}}^d=  \overline {Y}$ is a generalized minimal sufficient statistic for this family with respect to $\mathcal{L}_{\mathcal{J}}^{(\alpha)}$. 
Note that, for the function
\begin{equation*}
    h(t) = (-1)^{3t} \frac{3t}{2},
\end{equation*}
we have that 
$$\widetilde{E}_{\lambda}[h(T)] =\sum\limits_{y_1,y_2,y_3}(-1)^{3\bar y} \frac{3\bar y}{2} \widetilde{f}_{\lambda}(y_{1},y_2,y_3) = 0,$$ where $T=\overline {Y}$. Further, note that $\widetilde{P}_\lambda\{h(T)=0\}=\widetilde{P}_\lambda\{y_1,y_2,y_3:\sum_{i=1}^ny_i=0\}=\widetilde{f}_\lambda(0,0,0)= \frac{1-\lambda}{4}$ ($<1$). 
Thus, $T$ is not a generalized complete sufficient statistic for this family with respect to $\mathcal{L^{(\alpha)}_{J}}$.
\end{example}

\begin{remark}\label{2.R5.2} From the above example, one may note that a generalized minimal sufficient statistic is not necessarily generalized complete. Consequently, the converse part of Theorem \ref{2.T3.3} does not hold. $\hfill\Box$
\end{remark}
\hspace*{0.02 in}Although ${\frac{\overline{f}}{\overline{h}}}^d$ is not the complete sufficient statistic for $\mathcal{M}^{(\alpha)}$-family with respect to  $\mathcal{L^{(\alpha)}_{J}}$, one might be interested to know whether this is the generalized UMVUE for the given family with respect to $\mathcal{L^{(\alpha)}_{J}}$. We below give an example, which gives a negative answer to this question.
\begin{example}\label{2.Ex5.3}
     Consider the family of Bernoulli distributions as given in Example~\ref{2.Ex5.2}.
 Let $Y_1, Y_2$ be a random sample drawn from $f_\lambda$. Then, the deformed pdf associated with  $\mathcal{L^{(\alpha)}_{J}}$ is given by 
 \begin{equation*}
     \widetilde{f}_{\lambda}(y_1, y_2) = \frac{1-\lambda}{2}\Bigg[1+\overline{y}\frac{2\lambda-1}{1-\lambda}\Bigg],
 \end{equation*}
where $\overline{y}=\frac{y_1 + y_2}{2}$. Let $\widetilde{E}_\lambda[\overline{Y}] = \widetilde \tau(\lambda) $. Now, consider a function $\eta$ of $Y_1$ and $ Y_2$ as
\begin{equation*}
    \eta(Y_1, Y_2) = \begin{cases}
 1 , &\text{if $(Y_1, Y_2)= (0,0)$ or $(1,1)$}\\
-1, &\text{if $(Y_1, Y_2)= (0,1)$ or $(1,0)$}.
    \end{cases}
\end{equation*}
Then, $\widetilde{E}_\lambda[\eta(Y_1, Y_2)] = 0$, $\text{for all } \lambda$. Again, consider 
$$\widetilde{E}_\lambda[\eta(Y_1, Y_2)\cdot \overline{Y}] = \sum_{(y_1, y_2)}\overline{y}\cdot \eta(y_1, y_2)\cdot\widetilde{f}_\lambda(y_1, y_2) = \frac{1}{2}\left(\lambda - \frac{1}{2}\right),$$
which is not equal to zero except for $\lambda = 1/2$. Thus, by Theorem~\ref{2.T5.6}, we conclude that $\overline{Y}$ is not the generalized UMVUE for its expected value.
\end{example}
\section{Comparison between MDPDE and generalized UMVUE for $\mathcal{B^{(\alpha)}}$-family}\label{arna0}
As a robust alternative to classical estimation methods, DPD-based estimation methods have recently gained significant popularity among researchers. In the preceding section, we derive a methodology to obtain the generalized UMVUE for a given family of probability distributions with respect to a generalized likelihood function. On the other hand, there are several generalized estimators that are obtained based on some generalized likelihood functions. A natural question arises whether generalized UMVUE outperforms generalized estimators. The notion of deficiency, introduced by Hodges and Lehmann~\cite{HL}, can be used to provide a meaningful answer to this question. Hwang and Hu~\cite{HH} investigated the asymptotic expected deficiency (AED) of the maximum likelihood estimator (MLE)  relative to the UMVUE for a one-dimensional estimable function in an exponential family. In this study, we focus on $\mathcal{B^{(\alpha)}}$-family and derive the AED of the MDPDE relative to the generalized UMVUE. 
We begin this section by recalling the definition of the $\mathcal{B^{(\alpha)}}$-family, allowing readers to continue independently.

Consider a random sample $Y_1,\dots, Y_n$ drawn from a population with the probability density function $f_\lambda$ that belongs to the one-parameter $ \mathcal{B^{(\alpha)}}$-family, i.e., 
\begin{equation*}
  f_\lambda(y)=\begin{cases}
    \left[h(y)+Z(\lambda)+ {w(\lambda)}f(y)\right]^\frac{1}{\alpha-1}, & \text{for $y \in S$, $\lambda\in \Lambda,$}\\
    0, & \text{otherwise},
  \end{cases}
\end{equation*}
where $S$ is independent of the parameter and $f$ is a non-constant function. Further, we consider the generalized likelihood function 
$\mathcal{\ell}_{\mathcal{B}}^{(\alpha)}$ given by
\begin{equation}\label{2.E53}
    \mathcal{\ell}_{\mathcal{B}}^{(\alpha)}(y_1 ^n;\lambda):=\sum_{j=1} ^ {n} \left[\frac{ \alpha f_\lambda^{\alpha-1}(y_j) - 1}{\alpha-1}\right] - n \int f_\lambda^\alpha(y)  dy. 
\end{equation}
From the estimating equation given in \eqref{2.E9}, one may note that the MDPDE for $\mathcal{B^{(\alpha)}}$-family can be obtained by maximizing $\mathcal{\ell}_{\mathcal{B}}^{(\alpha)}$. Throughout this section, we assume $\mathcal{\ell}_{\mathcal{B}}^{(\alpha)}$ as the underlying generalized likelihood function.
\subsection{Asymptotic expected deficiency}\label{seSTH}
In this subsection, we derive the AED of MDPDE relative to generalized UMVUE, for a given one-parameter generalized estimable function of the $\mathcal{B^{(\alpha)}}$-family, and then compare the performance of MDPDE and generalized UMVUE.  By proceeding in the same lines as in Theorem \ref{2.T5.4}, we note that $T=\frac{1}{n}\sum_{j=1}^nf(Y_j)$ is a generalized complete sufficient statistic for $\mathcal{B^{(\alpha)}}$-family  
provided that the range of $w(\lambda)$ contains one-dimensional rectangle. Note that the MDPDE is a function of generalized sufficient statistic and has the invariance property (The proof follows in the same lines as in Zena~\cite{Zehna}). Thus, without any loss of generality, let $T$ be the MDPDE of $\lambda$. Consequently, for a generalized estimable function $\widetilde\tau$, $\widetilde\tau(T)$ is the MDPDE for $\widetilde\tau(\lambda)$.
Again, by Theorem \ref{2.T3.1}, we note that the generalized UMVUE of $\widetilde\tau(\lambda)$ is a function of $T$. Let $\widetilde U(T)$ be the generalized UMVUE of $\widetilde\tau(\lambda)$. 
If we consider the generalized risk function $\widetilde{R}_n$, under the squared error loss, as the measure of performance, then 
$$\widetilde{R}_n(\widetilde{\tau}(\lambda),  \widetilde{U}(T))=\widetilde{E}_\lambda[\widetilde{\tau}(\lambda)-\widetilde{U}(T)]^2$$
and
$$\widetilde{R}_n(\widetilde{\tau}(\lambda),  \widetilde{\tau}(T))=\widetilde{E}_\lambda[\widetilde{\tau}(\lambda)-\widetilde{\tau}(T)]^2.$$
Note that a statistical method, based on $n$ observations, is considered to be better compared to another method if it achieves the same risk as the latter has, but requires a larger sample size than $n$.
Let $k_n \;(\geq n)$ be the number of observations such that $\widetilde{R}_{n}(\widetilde{\tau}(\lambda),  \widetilde{U}(T))= \widetilde{R}_{k_n}(\widetilde{\tau}(\lambda),   \widetilde{\tau}(T))$, i.e., $\widetilde{\tau}(T)$ requires $k_n$ number of observations to achieve the same performance as that of $\widetilde{U}(T)$ with $n$ observations. Thus, the performance of one estimator relative to another estimator can be determined by considering either  $\lim_{n\rightarrow\infty} k_n/n$ (asymptotic relative efficiency (ARE)) or $\lim_{n\rightarrow\infty} (k_n-n)$ (asymptotic expected deficiency (AED)). Note that the  ARE is used most frequently due to its stability for a large sample. However, it fails when $\lim_{n\rightarrow\infty} k_n/n \rightarrow 1$. In that case, the AED can be used. However, for a given problem, obtaining AED is not that easy compared to ARE. 

It is reasonable to assume that the risks of $\widetilde{U}(T)$ and $\widetilde{\tau}(T)$ are strictly decreasing sequences in sample size with  $\widetilde{R}_n(\widetilde{\tau}(\lambda),  \widetilde{U}(T))\to 0$ and $ \widetilde{R}_n(\widetilde{\tau}(\lambda),  \widetilde{\tau}(T)) \to 0$, as $n\to \infty$. Thus, in line with Hwang and Hu~\cite{HH}, we assume that the generalized risks of $\widetilde{U}(T)$ and $\widetilde{\tau}(T)$ are of the form
\begin{equation}\label{2.E54}
 \widetilde{R}_n(\widetilde{\tau}(\lambda),  \widetilde{U}(T))= \frac{a}{n^r} + \frac{d}{n^{r+s}} + o(n^{-(r+s)}).  
\end{equation}
and
\begin{equation}\label{2.E55}
\widetilde{R}_n(\widetilde{\tau}(\lambda),  \widetilde{\tau}(T)) = \frac{a}{n^r} + \frac{b}{n^{r+s}} + o(n^{-(r+s)}),  
\end{equation}
respectively, where $r,s>0$ and $a,b,d \in \R$. Then, for $n\rightarrow\infty$, we have $\lim_{n\rightarrow\infty} k_n/n = 1$. On the other hand, the AED of $\widetilde{\tau}(T)$ relative to $\widetilde{U}(T)$, denoted by $AED[\widetilde{\tau}(T), \widetilde{U}(T)]$, is given by
\begin{eqnarray}\label{HTH}
AED[\widetilde{\tau}(T), \widetilde{U}(T)]=\lim_{n\rightarrow\infty} (k_n - n) = \begin{cases}
  \frac{b-d}{ra} & \text{if} \quad  s=1 \\
  \infty & \text{if} \quad  0<s<1 \\
  0 & \text{if} \quad  s>1.
\end{cases} 
\end{eqnarray}

In what follows, we derive the AED of MDPDE relative to generalized UMVUE for $\mathcal{B^{(\alpha)}}$-family.
To obtain the main result, we consider the $\mathcal{B^{(\alpha)}}$-family  with the following regularity conditions.
\begin{itemize}
    \item[$(a)$] $t$ is an interior point of $\Lambda$;
    \item[$(b)$] $\frac{d}{d\lambda}(\log N(\lambda))+n\lambda\frac{d}{d\lambda}(w^*(\lambda))=0$ and $\frac{d}{d\lambda}w^*(\lambda)>0$, for all $\lambda \in \Lambda$, where  $w^*(\lambda)=\frac{\alpha}{\alpha-1}w(\lambda)$ and $N(\lambda)=\bigintss \exp\left\{\frac{\alpha}{\alpha-1}\sum_{j=1}^nh(y_j)+w^*(\lambda)\sum_{j=1}^nf(y_j)\right\}dy_1^n$;
    \item[$(c)$] The density of $T$ can be differentiated with respect to $\lambda$ under the integral with respect to $t$ any number of times;
    \item[$(d)$] The functions $\widetilde{U}(t)$ and $\widetilde{\tau}(t)$ admit a convergent Taylor series for all the interior points of $\Lambda$.
\end{itemize}

Before stating the main result of this subsection, we below give a lemma wherein we drive the central moments of $T$.
\begin{lemma}\label{2.L5.1}
    For $i=1,2,\dots$, let $\widetilde{u}_i =\widetilde{E}_\lambda[T-\lambda]^i $ be the $i$-th central moment of $T$, and let $\widetilde{u}_0=\widetilde{E}_\lambda[T-\lambda]^0=1$. Then, under the regularity conditions (b) and (c), we have 
    $$\frac{d\widetilde u_i}{d\lambda} =  -i\widetilde u_{i-1}+n\left(\frac{d}{d\lambda}w^*(\lambda)\right)\widetilde u_{i+1}, \quad i=1,2,\dots. $$ Moreover, 
    $$\widetilde{u}_{2i-1}= O(n^{-i}) \quad \text{and} \quad \widetilde{u}_{2i}= O(n^{-i}), \quad  i=1,2, \dots .$$
\end{lemma}
\noindent\textbf{Proof :} Observe that $\widetilde{E}_\lambda[T-\lambda]^0=1$. Let $\widetilde{l}(y_1^n,\lambda)= \log \widetilde{f}_\lambda(y_1^n)$, where $\widetilde{f}_\lambda(y_1^n)$ is the underlying deformed probability density function.
Note that $\widetilde{E}_\lambda\left[\frac{d}{d\lambda}\widetilde{l}(y_1^n, \lambda)\right]=0$. Further, this along with the regularity condition $(b)$ implies that $\widetilde{E}_\lambda[T-\lambda]=0$.
Now, for $i=1,2,\dots$,
\begin{eqnarray*}
    &\frac{d \widetilde u_i}{d\lambda}& = - \int i (t-\lambda)^{i-1}N(\lambda)\exp \{k^*(t)+nw^*(\lambda)t\}d\mu(t)\\
    &&~~+\int \left\{{\frac{d}{d\lambda}(\log N(\lambda))}+nt\frac{d}{d\lambda}w^*(\lambda)\right\} (t-\lambda)^i N(\lambda)\exp \{k^*(t)+nw^*(\lambda)t\}d\mu(t),
\end{eqnarray*}
where $N(\lambda)\exp\{k^*(t)+nw^*(\lambda)t\}$ is the deformed probability density of $T$ with respect to some measure $\mu$. By using the regularity condition $(b)$, the above can equivalently be written as
\begin{eqnarray*}
    \frac{d\widetilde u_i}{d\lambda}
    =  -i\widetilde u_{i-1}+n\left(\frac{d}{d\lambda}w^*(\lambda)\right)\widetilde u_{i+1}, \quad i=1,2,\dots.
\end{eqnarray*}
Further, from the above equation, we get
\begin{equation}\label{2.E56}
  \widetilde u_2=\frac{1}{n \left(\frac{d}{d\lambda}w^*(\lambda)\right)},  
\end{equation}
\begin{equation}\label{2.E57}
\widetilde u_3= -\frac{\frac{d^2}{d\lambda^2}w^*(\lambda)}{n^2 \left(\frac{d}{d\lambda}w^*(\lambda)\right)^3 } ,  
\end{equation}
\begin{equation}\label{2.E58}
\widetilde u_4 = \frac{3}{\left(n\frac{d}{d\lambda}w^*(\lambda)\right)^2}+\frac{\left\{3\left(\frac{\frac{d^2}{d\lambda^2}w^*(\lambda)}{\frac{d}{d\lambda}w^*(\lambda)}\right)^2 - \frac{\frac{d^3}{d\lambda^3}w^*(\lambda)}{\frac{d}{d\lambda}w^*(\lambda)}\right\}}{\left(n\frac{d}{d\lambda}w^*(\lambda)\right)^3},
\end{equation}
and in general, 
$$~~~~~~~~~~~~~~~~~~~~~~~~~~~~~~~~~~~~~~\widetilde{u}_{2i-1}= O(n^{-i}) \quad \text{and} \quad \widetilde{u}_{2i}= O(n^{-i}), \quad  i=1,2, \dots .~~~~~~~~~~~~~~~~~~~~~~~~~~~~~~~~~~\hfill\Box$$

In the following theorem, we derive the AED of MDPDE ($\widetilde{\tau}(T)$) relative to generalized UMVUE ($\widetilde{U}(T)$) for $\widetilde{\tau}(\lambda)$.
\begin{theorem}\label{2.T6.5}
    Under the regularity conditions (a) - (d), the AED of $\widetilde{\tau}(T)$ relative to $\widetilde{U}(T)$, for $\widetilde{\tau}(\lambda)$, is given by
    \begin{eqnarray}
        AED[\widetilde{\tau}(T), \widetilde{U}(T)] = \frac{1}{\frac{d}{d\lambda}w^{*}(\lambda)}\left\{\frac{\frac{d^3}{d\lambda^3}\widetilde{\tau}(\lambda)}{\frac{d}{d\lambda}\widetilde{\tau}(\lambda)}+\frac{1}{4}\left(\frac{\frac{d^2}{d\lambda^2}\widetilde{\tau}(\lambda)}{\frac{d}{d\lambda}\widetilde{\tau}(\lambda)}\right)^2\right\} - \frac{\frac{d^2}{d\lambda^2}w^{*}(\lambda)}{\left(\frac{d}{d\lambda}w^{*}(\lambda)\right)^2} \frac{\frac{d^2}{d\lambda^2}\widetilde{\tau}(\lambda)}{\frac{d}{d\lambda}\widetilde{\tau}(\lambda)}.
    \end{eqnarray}
\end{theorem}
\noindent \textbf{Proof:} 
Using the regularity condition $(d)$, we can write
\begin{equation}\label{2.E59}
    \widetilde{U}(T) = \widetilde{U}(\lambda)+ \sum_{i=1}^\infty \left(\frac{d^i}{d\lambda^i}\widetilde{U}(\lambda)\right)\frac{(T-\lambda)^i}{i!}
\end{equation}
and
\begin{equation}\label{2.E60}
 \widetilde{\tau}(T) = \widetilde{\tau}(\lambda)+ \sum_{i=1}^\infty \left(\frac{d^i}{d\lambda^i}\widetilde{\tau}(\lambda)\right)\frac{(T-\lambda)^i}{i!}.  
\end{equation}
Again, the generalized risk of the estimator $\widetilde{U}(T)$ can be written as
\begin{equation}\label{2.E61}
\widetilde{E}_{\lambda}[\widetilde{U}(T)-\widetilde{\tau}(\lambda)]^2= \widetilde{E}_{\lambda}[\widetilde{U}(T)-\widetilde{U}(\lambda)]^2-[\widetilde{U}(\lambda)-\widetilde{\tau}(\lambda)]^2.
\end{equation}
Now, from \eqref{2.E59} and Lemma \ref{2.L5.1}, we get
\begin{eqnarray}\label{2.E62}
    &\widetilde{E}_\lambda[\widetilde{U}(T)-\widetilde{U}(\lambda)]^2 &= \widetilde{E}_\lambda\left[\sum_{i=1}^\infty \left(\frac{d^i}{d\lambda^i}\widetilde{U}(\lambda)\right)\frac{(T-\lambda)^i}{i!}\right]^2 \quad \nonumber \\
    &&= \left(\frac{d}{d\lambda}\widetilde{U}(\lambda)\right)^2 \widetilde{u}_2 + \left(\frac{d}{d\lambda}\widetilde{U}(\lambda)\right)\left(\frac{d^2}{d\lambda^2}\widetilde{U}(\lambda)\right)\widetilde{u}_3 \nonumber\\
    && + \left\{\frac{1}{4}\left(\frac{d^2}{d\lambda^2}\widetilde{U}(\lambda)\right)^2 +\frac{1}{3}\left(\frac{d}{d\lambda}\widetilde{U}(\lambda)\right)\left(\frac{d^3}{d\lambda^3}\widetilde{U}(\lambda)\right)\right\}\widetilde{u}_4 + O(n^{-3})
\end{eqnarray}
and
\begin{eqnarray}\label{2.E63}
&[\widetilde{U}(\lambda)-\widetilde{\tau}(\lambda)]^2& = \left\{\widetilde{E}_\lambda\left[\sum_{i=1}^\infty \left(\frac{d^i}{d\lambda^i}\widetilde{U}(\lambda)\right)\frac{(T-\lambda)^i}{i!}\right]\right\}^2 \nonumber \\
&& =  \frac{1}{4}\left(\frac{d^2}{d\lambda^2}\widetilde{U}(\lambda)\right)^2\widetilde{u}_2^2+O(n^{-3}).
\end{eqnarray}
Further, by taking expectation on both sides of \eqref{2.E59} and then differentiating with respect to $\lambda$, we get
\begin{equation*}
    \frac{d}{d\lambda}\widetilde{\tau}(\lambda) = \sum_{i=1}^\infty \frac{n\left(\frac{d}{d\lambda}w^*(\lambda)\right)\left(\frac{d^i}{d\lambda^i}\widetilde{U}(\lambda)\right)\widetilde{u}_{i+1}}{i!}. 
\end{equation*}
Upon using Lemma \ref{2.L5.1} in the above equation, we get
\begin{equation}\label{2.E64}
    \frac{d \widetilde{U}(\lambda)}{d\lambda} = \frac{d \widetilde{\tau}(\lambda) }{d\lambda} + \frac{1}{2n\left(\frac{d}{d\lambda}w^*(\lambda)\right)}{\left\{-\left(\frac{d^3}{d\lambda^3}\widetilde{U}(\lambda)\right) +\frac{ \left(\frac{d^2}{d\lambda^2}\widetilde{U}(\lambda)\right)\left(\frac{d^2}{d\lambda^2}w^*(\lambda)\right)}{\left(\frac{d}{d\lambda}w^*(\lambda)\right)}\right\}}+O(n^{-2}),
\end{equation}
\begin{equation}\label{2.E65}
 \frac{d^2}{d\lambda^2}\widetilde{U}(\lambda) = \frac{d^2}{d\lambda^2}\widetilde{\tau}(\lambda) + O(n^{-1}), 
\end{equation}
and 
\begin{equation}\label{2.E66}
 \frac{d^3}{d\lambda^3}\widetilde{U}(\lambda) = \frac{d^3}{d\lambda^3}\widetilde{\tau}(\lambda) + O(n^{-1}).     
\end{equation}
Then, by using \eqref{2.E64}, \eqref{2.E65} and \eqref{2.E66} in \eqref{2.E62} and \eqref{2.E63}, we get,  from \eqref{2.E61},
\begin{eqnarray}\label{2.E67}
\widetilde{E}_{\lambda}[\widetilde{U}(T)-\widetilde{\tau}(\lambda)]^2&=& \frac{\left(\frac{d}{d\lambda}\widetilde{\tau}(\lambda)\right)^2}{n\left(\frac{d}{d\lambda}w^*(\lambda)\right)}+\frac{1}{2} \frac{\left(\frac{d^2}{d\lambda^2}\widetilde{\tau}(\lambda)\right)^2}{\left(n\frac{d}{d\lambda}w^*(\lambda)\right)^2}+O(n^{-3})\nonumber\\&=&\frac{a}{n}+\frac{d}{n^2}+O(n^{-3}),
\end{eqnarray}
where \begin{eqnarray}\label{x2}
    a= \frac{\left(\frac{d}{d\lambda}\widetilde{\tau}(\lambda)\right)^2}{\left(\frac{d}{d\lambda}w^*(\lambda)\right)}\text{ and }
d= \frac{1}{2} \frac{\left(\frac{d^2}{d\lambda^2}\widetilde{\tau}(\lambda)\right)^2}{\left(\frac{d}{d\lambda}w^*(\lambda)\right)^2}.
\end{eqnarray}
Further, in view of \eqref{2.E60} and Lemma \ref{2.L5.1}, the generalized risk of $\tilde{\tau}(T)$ is given by
\begin{eqnarray}\label{x34}
\widetilde{E}_\lambda[\widetilde{\tau}(T)-\widetilde{\tau}(\lambda)]^2 & =& \widetilde{Var}_\lambda[\widetilde{\tau}(T)-\widetilde{\tau}(\lambda)] + \{\widetilde{E}_\lambda[\widetilde{\tau}(T)-\widetilde{\tau}(\lambda)]\}^2 \nonumber \\
    & = &\widetilde{Var}\left[\sum_{i=1}^\infty \left(\frac{d^i}{d\lambda^i}\widetilde{\tau}(\lambda)\right)\frac{(T-\lambda)^i}{i!}\right] + \left\{\widetilde{E}_\lambda\left[\sum_{i=1}^\infty \left(\frac{d^i}{d\lambda^i}\widetilde{\tau}(\lambda)\right)\frac{(T-\lambda)^i}{i!}\right]\right\}^2 \nonumber \\
    & = &\frac{\left(\frac{d}{d\lambda}\widetilde{\tau}(\lambda)\right)^2}{n\left(\frac{d}{d\lambda}w^*(\lambda)\right)} + \Bigg\{\left(\frac{d}{d\lambda}\widetilde{\tau}(\lambda)\right)\left(\frac{d^3}{d\lambda^3}\widetilde{\tau}(\lambda)\right) - \left(\frac{d}{d\lambda}\widetilde{\tau}(\lambda)\right)\left(\frac{d^2}{d\lambda^2}\widetilde{\tau}(\lambda)\right)\frac{\left(\frac{d^2}{d\lambda^2}w^*(\lambda)\right)}{\left(\frac{d}{d\lambda}w^*(\lambda)\right)} \nonumber \\
    && + \frac{3}{4}\left(\frac{d^2}{d\lambda^2}\widetilde{\tau}(\lambda)\right)^2 \Bigg\} \frac{1}{n^2\left(\frac{d}{d\lambda}w^*(\lambda)\right)^2}+O(n^{-3})\nonumber
    \\&=&\frac{a}{n}+\frac{b}{n^2}+O(n^{-3}),
\end{eqnarray}
where $a$ is the same as in \eqref{x2} and 
\begin{eqnarray*}
b&=& \Bigg\{\left(\frac{d}{d\lambda}\widetilde{\tau}(\lambda)\right)\left(\frac{d^3}{d\lambda^3}\widetilde{\tau}(\lambda)\right) - \left(\frac{d}{d\lambda}\widetilde{\tau}(\lambda)\right)\left(\frac{d^2}{d\lambda^2}\widetilde{\tau}(\lambda)\right)\frac{\left(\frac{d^2}{d\lambda^2}w^*(\lambda)\right)}{\left(\frac{d}{d\lambda}w^*(\lambda)\right)} \nonumber 
     + \frac{3}{4}\left(\frac{d^2}{d\lambda^2}\widetilde{\tau}(\lambda)\right)^2 \Bigg\} \frac{1}{\left(\frac{d}{d\lambda}w^*(\lambda)\right)^2}.
    \end{eqnarray*}
    Note that \eqref{2.E67} and \eqref{x34} are in the forms of \eqref{2.E54} and \eqref{2.E55}, respectively. Finally, by putting the values of $a, b$ and $d$ in \eqref{HTH}, we get the required result.
\subsection{ An application in stress-strength reliability}
In reliability theory, the stress-strength model is one of the commonly used frameworks. It is used to evaluate the system efficiency. The stress-strength reliability measures the probability that the strength ($Y$)  of a system exceeds the imposed stress ($X$) on it. A system fails when stress exceeds its strength. The stress-strength reliability model, introduced by Birnbaum~\cite{B2}, is widely used in different fields, including reliability engineering, statistics, and clinical trials. For example, Hauck et al.~\cite{HHA} used stress-strength reliability as a measure of the effect of treatment, where $Y$ represents the response of a control group to a therapeutic approach and $X$ represents the reaction of a treated group. 
The estimation problem of stress-strength reliability has been studied for several classical distributions, including exponential, gamma, and normal; see, for example, Church and Harris~\cite{CH}, Huang et al.~\cite{HMW}, Kundu and Gupta~\cite{KG}, and the references therein. 
We below discuss a generalized estimation problem for stress-strength model, where both stress and strength have the Student distribution.

Let $Y$ be a random variable representing the strength of a system with the probability density function
given by 
\begin{equation}\label{1bn}
  f_\lambda(y) = \frac{\Gamma(\frac{\nu+1}{2})}{\Gamma(\frac{\nu}{2})\sqrt{\pi \nu }} \left[1+\frac{1}{\nu} ({y-\mu})^2\right]^{\frac{1}{\alpha-1}},\quad y\in \mathbb{R},
  \end{equation}
  where $\nu~(>2)$  is fixed, $\sigma^2=1$ and $\mu$ is unknown.
 Note that $  f_\lambda \in\mathcal{B^\alpha}$- family with $w^*(\mu)=-\frac{2\alpha\mu}{\nu(\alpha-1)}\left(\frac{\Gamma\frac{\nu+1}{2}}{\Gamma \frac{\nu}{2}\sqrt{\pi \nu}}\right)^{\alpha-1}$. Further, let $X$ be another random variable representing the stress imposed on the same system with the probability density function given by
 \begin{equation}\label{2bn}
  h_\lambda(x) = \frac{\Gamma(\frac{\nu+1}{2})}{\Gamma(\frac{\nu}{2})\sqrt{\pi \nu }} \left[1+\frac{x^2}{\nu} \right]^{\frac{1}{\alpha-1}}, \quad x\in \mathbb{R},
  \end{equation}
  where $\nu$ is the same as in \eqref{1bn}. Let $X_1,\dots,X_n$ and $Y_1,\dots,Y_n$ be independent random samples drawn from the distributions given in \eqref{1bn} and \eqref{2bn}, respectively. Then, the deformed probability distributions corresponding to \eqref{1bn} and \eqref{2bn} are given by 
  \begin{equation}
    \widetilde{f}_\lambda(y_1^n)= \frac{\exp\left\{-\frac{1}{2}\sum_{j=1}^{n}\left(\frac{y_j-\mu}{\sigma^*}\right)^2\right\}}{\bigints\exp\left\{-\frac{1}{2}\sum_{j=1}^{n}\left(\frac{y_j-\mu}{\sigma^*}\right)^2\right\}dy_1^n} 
  \end{equation}
  and 
  \begin{equation}
    \widetilde{h}_\lambda(x_1^n)= \frac{\exp\left\{-\frac{1}{2}\sum_{j=1}^{n}\left(\frac{x_j}{\sigma^*}\right)^2\right\}}{\bigints\exp\left\{-\frac{1}{2}\sum_{j=1}^{n}\left(\frac{x_j}{\sigma^*}\right)^2\right\}dx_1^n}, 
  \end{equation}
respectively, where $\sigma^*=\left[\frac{2\alpha}{1-\alpha }\frac{1}{{(\sqrt{\nu})}^{1+\alpha}}{\left(\frac{\Gamma{\frac{\nu+1}{2}}}{\Gamma{\frac{\nu}{2}}\sqrt{\pi}}\right)}^{\alpha-1}\right]^{-1/2}$. Consequently, the deformed marginal probability distributions of $Y_1$ and $X_1$  are given by
$$\widetilde{f}_\lambda(y_1)=\frac{\exp\left\{-\frac{1}{2}\left(\frac{y_1-\mu}{\sigma^*}\right)^2\right\}}{\bigints\exp\left\{-\frac{1}{2}\left(\frac{y_1-\mu}{\sigma^*}\right)^2\right\}dy_1}$$
and 
$$\widetilde{h}_\lambda(x_1)=\frac{\exp\left\{-\frac{1}{2}\left(\frac{x_1}{\sigma^*}\right)^2\right\}}{\bigints\exp\left\{-\frac{1}{2}\left(\frac{x_1}{\sigma^*}\right)^2\right\}dx_1},$$
respectively.
Then, $\widetilde{P}_\lambda\{Y_1<X_1\}=\Phi\left(\frac{\mu}{\sqrt{2}\sigma^*}\right)$, which is a function of $\mu$; here $\Phi$ is the cumulative distribution function of the standard normal distribution. By Theorem 11 of Gayen and Kumar \cite{Projection_theorems}, we have that the MDPDE of $\mu$ is $\overline{Y}$. Consequently, by using the invariance property, we get that $\Phi\left(\frac{\overline{Y}}{\sqrt{2}\sigma^*}\right)$ is the MDPDE  of $\Phi\left(\frac{\mu}{\sqrt{2}\sigma^*}\right)$. Further, note that the range of $w^*(\mu)$ contains one-dimentional rectangle. Thus, in line with Theorem \ref{2.T5.4}, we can show that $\overline{Y}=\frac{1}{n}\sum_{j=1}^n Y_j$ is a generalized complete sufficient statistic for $\mu$. Note that $\widetilde{E}_\lambda\left[\Phi\left(\frac{\sqrt{n}}{\sqrt{2n-1}}\frac{\overline{Y}}{\sigma^*}\right)\right]= \Phi\left(\frac{\mu}{\sqrt{2}\sigma^*}\right)$. This implies that $\Phi\left(\frac{\sqrt{n}}{\sqrt{2n-1}}\frac{\overline{Y}}{\sigma^*}\right)$ is a generalized unbiased estimator which is a function of $\overline Y$. Consequently, by Theorem~\ref{2.T3.1}, we get that $\Phi\left(\frac{\sqrt{n}}{\sqrt{2n-1}}\frac{\overline{Y}}{\sigma^*}\right)$  is the generalized UMVUE for $\Phi\left(\frac{\mu}{\sqrt{2}\sigma^*}\right)$. Now, by putting $\widetilde{\tau}(\lambda)=\Phi\left(\frac{\mu}{\sqrt{2}\sigma^*}\right)$ and $w^*(\lambda)=w^*(\mu)$ in Theorem \ref{2.T6.5}, the AED of MDPDE relative to generalized UMVUE, for $\Phi\left(\frac{\mu}{\sqrt{2}\sigma^*}\right)$, is given by
\begin{equation*}
    AED\left[\Phi\left(\frac{\overline{Y}}{\sqrt{2}\sigma^*}\right), \Phi\left(\frac{\sqrt{n}}{\sqrt{2n-1}}\frac{\overline{Y}}{\sigma^*}\right)\right] = \frac{(4+\sigma^*)\mu^2-8{\sigma^*}^2}{16\sigma^*}.
\end{equation*}
From the above expression, we conclude the following.
\begin{enumerate}
    \item [$(i)$] If $|\mu|<\sqrt{\frac{8}{4+\sigma^*}}\sigma^*$, then the MDPDE is better than the generalized UMVUE for $\Phi\left(\frac{\mu}{\sqrt{2}\sigma^*}\right)$;
    \item[$(ii)$] If $|\mu|>\sqrt{\frac{8}{4+\sigma^*}}\sigma^*$, then the generalized UMVUE is better than  the MDPDE for $\Phi\left(\frac{\mu}{\sqrt{2}\sigma^*}\right)$.
\end{enumerate}
Further, note that, if the degrees of freedom of the underlying Student distributions varies from $3$ to $30$, then $\sigma^*$ varies from $0.95$ to $1$, and consequently, $0.19<\Phi\left(\frac{\mu}{\sqrt{2}\sigma^*}\right)<0.81$ when $|\mu|<\sqrt{\frac{8}{4+\sigma^*}}\sigma^*$, and $\Phi\left(\frac{\mu}{\sqrt{2}\sigma^*}\right)<0.19$ or $\Phi\left(\frac{\mu}{\sqrt{2}\sigma^*}\right)>0.81$ when $|\mu|>\sqrt{\frac{8}{4+\sigma^*}}\sigma^*$. In a given stress-strength model, the underlying system functions as long as the strength of the system is greater than the stress imposed on it. In this case, the stress-strength reliability, $\Phi\left(\frac{\mu}{\sqrt{2}\sigma^*}\right)$, is close to $1$, which implies that the generalized UMVUE is better than the MDPDE. 
\section{Conclusion}\label{2.S6}
The notion of sufficient statistics is majorly useful in data reduction. It helps in finding a better estimator by make use of Rao-Blackwell theorem. The notion of sufficiency was initially defined for the ML estimation. However, one of the inherent drawbacks of the ML estimation is that it gets affected by outliers and fails to give desired estimate. In order to reduce the effect of outliers, several divergence-based estimation methods (namely, MDPDE, MLDPDE, etc.) were developed in the literature. These methods contain one or more tuning parameter(s) that balances the robustness and the efficiency of an estimator. These estimation methods are associated with some generalized likelihood functions. Based on these likelihood functions, the notion of generalized sufficiency was recently introduced in the literature. Note that the sufficiency alone does not provide the best estimator. Thus, in line with the generalized sufficiency, we introduce the notion of generalized complete statistic and generalized ancillary statistic for a generalized likelihood estimation. Subsequently, we derive the generalized Lehmann-Scheffé theorem that is useful to find the generalized UMVUE for a given family of distributions. Further, we obtain the generalized complete statistic for $\mathcal{B^{(\alpha)}}$-family and subsequently, derive the generalized UMVUE for this family. We also show that the underlying generalized minimal sufficient statistic for $\mathcal{M^{(\alpha)}}$-family is not complete. 
Further,  we derive the AED of MDPDE relative to generalized UMVUE, for $\mathcal{B^{(\alpha)}}$-family, and compare the performance of MDPDE and generalized UMVUE. Finally, we give an application of the developed results in stress-strength reliability model.

Here we study some classical inference problems for a generalized likelihood estimation. The study of some Bayesian inference problems under the setup of a generalized likelihood estimation is yet to be explored. We are currently working on this problem. If we succeed, the findings will be reported in a future communication.
 \\\\{\bf Acknowledgments}\\
\hspace*{0.2 in}The first author sincerely acknowledges the financial support received from the Ministry of Education, Government of India, under the Prime Minister's Research Fellowship (PMRF) scheme. 
\\\\{\bf Conflicts of Interest}\\
\hspace*{0.2 in} The authors declare no conflict of interest.
 
\end{document}